\documentclass{article}

\usepackage{arxiv}
\usepackage[utf8]{inputenc} 
\usepackage[T1]{fontenc}    
\usepackage{hyperref}       
\usepackage{url}            
\usepackage{booktabs}       
\usepackage{amsfonts}       
\usepackage{amsmath, amssymb}
\usepackage{graphicx}
\usepackage{natbib}
\usepackage{doi}
\usepackage[nolist]{acronym}
\usepackage[capposition=top]{floatrow}
\usepackage{multirow}
\usepackage{threeparttable}
\usepackage[title]{appendix}%
\usepackage{mathtools}

\title{Design of Targeted Community-Based Resource Allocation in the Presence of Vaccine Hesitancy via a Data-Driven Compartmental Stochastic Optimization Model}

\author{Hieu Bui\\
        Department of Industrial Engineering\\
        University of Arkansas\\
        Fayetteville, AR 72701\\
        \And
        Sandra Ek\c{s}io\u{g}lu\thanks{Corresponding author}\\
        Department of Industrial Engineering\\
        University of Arkansas\\
        Fayetteville, AR 72701\\
        \And
        Rub\'en Proa\~{n}o\\
        Department of Industrial and Systems Engineering\\
        Rochester Institute of Technology\\
        Rochester, NY 14623\\
        \And
        Haoming Shen\\
        Department of Industrial Engineering\\
        University of Arkansas\\
        Fayetteville, AR 72701\\
        }

\renewcommand{\shorttitle}{Resource Allocation via Compartmental MSP}

\hypersetup{
pdftitle={Design of Targeted Community-Based Resource Allocation in the Presence of Vaccine Hesitancy via a Data-Driven Compartmental Stochastic Optimization Model},
pdfsubject={},
pdfauthor={Hieu Bui, Sandra Ek\c{s}io\u{g}lu, Rub\'en Proa\~{n}o, Haoming Shen},
pdfkeywords={compartmental model in epidemiology, compartmental multi-stage stochastic program, decision making under uncertainty, resource allocation, vaccine hesitancy},
}

\begin{document}

\maketitle

\begin{abstract}
    Vaccines have proven effective in mitigating the threat of severe infections and deaths during outbreaks of infectious diseases. However, the reluctance or refusal to get vaccinated, known as \ac{VH}, presents a significant challenge in predicting the spread of the disease and assessing the need for healthcare resources in different regions and population groups. We propose a modeling framework that integrates an epidemiological compartmental model that captures the spread of an infectious disease within a \acf{MSP} that determines the allocation of critical resources under uncertainty. The proposed compartmental MSP model adaptively manages the allocation of resources to account for changes in population behavior toward vaccines (i.e., variability in \ac{VH}), the unique patterns of disease spread, and the availability of healthcare resources over time and space. The compartmental MSP model allowed us to analyze the price of fairness in resource allocation. We developed a case study using real-life data about COVID-19 vaccination uptake from January to May 2021 and data about healthcare resources in Arkansas, U.S. Using our numerical analysis, we observed that ($i$) delaying the initial deployment of additional ventilators (critical resource) by one month could lead to an average increase in the expected number of deaths by 285.41/month, highlighting the importance of prompt action; ($ii$) each additional ventilator in the initial stockpile and in supply lead to a decrease in the expected number of deaths by 1.09/month and 0.962/month, respectively, highlighting the importance of maintaining a large stockpile and scalable production response; ($iii$) the price of ensuring equitable allocation of resources changes over time and space. This price is highest during the peak of a disease outbreak and in highly populated regions. In summary, this study highlights the importance of flexible and informed decision-making in public health crises, and the crucial role of advanced planning and preparedness for future outbreaks, offering a framework for effective resource allocation adaptable to various public health emergencies.
\end{abstract}


\section{Introduction} \label{s:intro}

Infectious disease outbreaks with high reproduction numbers ($R_O$), such as Measles ($R_O=16$), COVID-19 ($R_0\in [1.35, 2.11]$), etc., can affect a large portion of the population in a short period of time \citep{sy_2021}. Therefore, the implementation and effective execution of public health strategies, such as vaccination campaigns, mask mandates, and social distancing, can greatly reduce the spread of the disease and save lives \citep{moghadas_2021, brooks_2021, chu_2020}. However, the success of vaccination campaigns, which are crucial to controlling disease outbreaks, is often impeded by \acf{VH}. This reluctance or refusal to get vaccinated leads to a disparity in vaccination rates and prolongs the repercussions of the disease, as in the case of the COVID-19 pandemic. Furthermore, \ac{VH} complicates the prediction of vaccine demand and the allocation of critical resources \citep{blasioli_2023}.

Making well-informed and timely decisions is essential to effectively manage disease outbreaks. Compartmentalized epidemiological models, like the \acf{SIR} and \acf{SEIR}, play a critical role in forecasting the progression of infectious diseases, thus aiding in the design of containment strategies. The \ac{SIR} model divides the population into three groups: Susceptible, Infectious, and Recovered, offering a basic framework for the dynamics of the disease \citep{kermack_1927}. The \ac{SEIR} model expands on this by including an `Exposed' group, which is crucial for understanding the progression of diseases that have a latency period, such as COVID-19. These models provide insight into disease transmission patterns that are vital to planning and implementing public health responses \citep{kong_2022, xiang_2021}. The traditional \ac{SIR} and \ac{SEIR} models typically assume constant transition rates between compartments, such as infection and recovery rates, yet this simplification does not always hold true. Advanced models address this by incorporating the uncertainty of these model parameters. However, even these adaptations of the \ac{SEIR} model that utilize stochastic compartments to inform COVID-19 strategies often neglect critical aspects of healthcare capacity and \ac{VH} \citep{buonomo_2022, olivera_2023}. Several researchers use operations research (OR) models to determine resource allocation strategies, which are critical decisions for healthcare management during disease outbreaks. These models are used to optimize the distribution of limited resources, such as ventilators during the COVID-19 pandemic, ensuring effective utilization \citep{mehrotra_2020, bhavani_2021, koonin_2020, bertsimas_2021}. From the wide range of OR models, \ac{SP} and \acf{MSP} provide a robust framework for adaptive decision-making under uncertainty \citep{birge2011introduction}. 

Incorporating dynamic disease progression data from compartmental models would significantly enhance the effectiveness of both the SP and MSP models. However, the integration of these models has been impeded by the computational challenges associated with solving each of them. Only a few studies have attempted to integrate these models to enhance the quality of resource allocation decisions \citep{yin_2021, yin_2023}. By incorporating the predictive power of compartmental models into SP and MSP models, we can identify and evaluate resource allocation strategies that are informed by the specifics of disease spread and progression and the availability of healthcare resources. This realization leads to our first research question (RQ1): \textit{``Can we develop an effective resource allocation model of critical resources needed during disease outbreaks?''} To demonstrate how we can tackle this, we develop an MSP compartmental model for the COVID-19 pandemic. We first develop a compartmental model that is an extension of the \ac{SEIR} model to capture the specific dynamics of the spread of the COVID-19 virus, the impact of \ac{VH} on the number of people exposed, and the availability of healthcare resources on the expected number of people recovered and of deaths. The outcomes of the modified \ac{SEIR} model are incorporated in an \ac{MSP} model that determines resource allocation.

During disease outbreaks and other health emergencies, resource allocation becomes a critical concern. For example, the COVID-19 pandemic caused an urgent demand for essential medical equipment, such as ventilators for people in critical condition and personal protective equipment. Recognizing the gravity of the situation, the U.S. government invoked the \acf{DPA}, urging car manufacturers to pivot their production lines to manufacture ventilators \citep{albergotti_2020}. This proactive measure showcased the influence of government interventions in crisis management. However, it also prompts reflection on its execution and raises questions such as: Could different results have been achieved with an earlier or more intensive intervention? As we anticipate and prepare for potential future outbreaks, it becomes imperative to assess the timing and intensity of interventions on healthcare outcomes. This leads to the following research question (RQ2): \textit{``What should be effective resource allocation strategies at the start of an outbreak?''} Specifically, how does the timing and scope of interventions, such as \ac{DPA}, impact the spread of the disease? How does VH impact the effectiveness of resource allocation? By exploring this research question, our aim is to provide guidance on optimizing resource allocation strategies to ensure that they are timely and effective.

Ensuring access to healthcare resources is crucial, especially during disease outbreaks, since strategies used for resource allocation can have significant impacts on the affected population. Several OR models take a utilitarian approach by focusing on minimizing costs or maximizing the benefits from the allocation of limited healthcare resources. Such an approach allocates resources to highly populated areas and those that have access to healthcare care, which may lead to inequity. In the U.S., access to healthcare varies greatly, and rural and marginalized communities often face more challenges in accessing the necessary services and resources \citep{mishra_2021}. These disparities can affect the effectiveness and success of public health intervention strategies. This leads to our next research question (RQ3): \textit{``What are the trade-offs between equity, equality, and the effectiveness of resource allocation?''} Exploring this question will help to understand how resource allocation strategies influence disease outbreak management, both in the short and long term.

\section{Literature review} \label{s:lit}

We identified three streams of literature related to our research: compartmental models in epidemiology, stochastic optimization models for resource allocation, and modeling equity in healthcare. The first stream of literature focuses on compartmental models, such as the \ac{SIR} and \ac{SEIR}. These models categorize populations into different subgroups according to their health status. These models provide a structured approach to assess disease progression and assess the impact of interventions on the affected population. The second stream of literature focuses on multistage stochastic programs (MSPs), models that aid in making decisions under uncertainty. Our study integrates compartmental models within an \ac{MSP} framework to optimize resource allocation to manage and contain infectious diseases. The third stream addresses equity in resource allocation. The goal is to ensure a fair and accessible distribution of healthcare resources, which is vital during infectious disease outbreaks. Each of the following subsections provides a comprehensive overview of the literature. We highlight their relevance to our proposed research model.

\subsection{Compartmental Models in Epidemiology}

Epidemiological models are key to understanding and predicting the course of infectious disease outbreaks. With the emergence of COVID-19, the field of mathematical epidemiology quickly attracted broader interest.  Taking advantage of existing approaches, researchers adapted and applied various methods specifically to study COVID-19, including compartmental models, structured metapopulations, agent-based networks, deep learning, and complex networks \citep{wu_2020,adiga_2020,rodriguez_2021,zhou_2020,chang_2021,kerr_2021,chang_2021b}. Among these, compartmental models are essential in infectious disease research due to their effective representation of the dynamics of disease transmission. In these models, the population is divided into several subgroups or compartments, each representing a specific stage in the progression of the disease. Consider, for example, the SIR model that divides the population into three compartments: Susceptible, Infectious, and Recovered. This model provides a basic framework for tracking disease spread by monitoring changes in population proportions within these compartments. Building on the \ac{SIR} model, the \ac{SEIR} model adds an Exposed compartment, highlighting a critical phase where individuals have contracted the disease but are not yet infectious. This addition is particularly important for diseases that have an incubation period, such as COVID-19. Several researchers extended the \ac{SEIR} model to capture the unique transmission dynamics of the COVID-19 virus, as well as the impacts of vaccination and public health mandates on the spread of the disease. Kong et al. (\citeyear{kong_2022}) provide a detailed overview of these extensions, which consider various stages of disease progression, intervention strategies, and demographic factors \citep{kimathi_2021, tuite_2020, foy_2021}.

Classic epidemiological models assume that model parameters, such as the transmission rate, the recovery rate, and the contact rate, are constant throughout the course of the epidemic being modeled. Most recently, epidemiological models have been updated to incorporate the variability and uncertainty of these parameters. Models that take into account the randomness that exists in the system can lead to better decision-making. For example, \citep{lekone_2006,gray_2011,faranda_2020, ganyani_2021, maltsev_2021, gatto_2020, giordano_2020, kretzschmar_2020, kucharski_2020} propose stochastic epidemiological models. These models use historical data to estimate how parameters change over time, and thus, provide a more realistic representation of the dynamics of infectious diseases.

VH has consistently presented a significant challenge to the healthcare system, a concern that has been recognized from the \ac{WHO} as a global health threat (\citep{who_2019}. The advent of COVID-19 vaccines once again highlighted this issue. Several compartmental models capture the impact of vaccine availability and hesitancy on transmission rate. \cite{buonomo_2022} demonstrate that voluntary vaccination alone cannot stop the spread of the disease. \cite{choi_2020} develop a game-theoretic epidemiological model to explore vaccination and social distancing as strategies to control the spread of the disease. They use a group behavioral model to determine optimal strategies for individuals. \cite{awad_2021} propose a compartmental model that uses geospatial and \ac{VH} information to develop strategies to distribute vaccines. \cite{jing_2023} examine the role of \ac{VH} in the emergence of new SARS-CoV-2 variants. \cite{BER_24} present a compartmental model that captures the impact of \ac{VH} on disease outcomes, such as the expected number of hospitalizations and deaths. The model runs for different values of \ac{VH} (which are taken from a certain distribution) to estimate the impacts of uncertain \ac{VH} on model outcomes. The authors performed a sensitivity analysis to assess the impact of the availability of healthcare resources and vaccines on the expected number of deaths. Despite the variety of compartmental models available, none examines the impacts of the temporal and spatial variations of \ac{VH} on demand for healthcare resources. This gap underscores the need for models that use predictions of temporal and spatial disease spread to determine the effective allocation of critical resources.

Based on our review of the literature, epidemiological models are typically used to support decision making via ``What-if" analysis. This approach explores how changes in model parameters or assumptions affect the spread of the disease and the need for resources. Some researchers integrate epidemiological models into optimization models to support optimal decision making during disease outbreaks. For example, \cite{abdin_2023} propose a non-linear programming formulation of a compartmental model to present the complex dynamics of infectious diseases, offering a nuanced perspective on disease transmission. \cite{yin_2021} propose an MSP compartmental model to support resource allocation during the Ebola outbreak in West Africa. The MSP uses information on stochastic disease progression to determine resource allocation to control the infectious disease outbreak. The model determines the fair allocation of critical resources. 

Our study investigates the impact of vaccination and VH on disease transmission dynamics. The model we propose enables decision makers to assess the impact that the timing, scale, and scope of resource allocation have on the spread of the disease and disease outcomes. While previous research has addressed aspects of resource management during disease outbreaks, none has presented a comprehensive, data-driven model, as we propose. Our model highlights the advantages of merging the predictive capabilities of compartmental models with optimization techniques, providing valuable insights into effective disease control strategies.

\subsection{Stochastic Optimization Models for Resource Allocation}

MSP models are used to support sequential decision making under uncertainty. The simplest MSP model, the two-stage stochastic program (2-SP), is commonly used to minimize expected costs (maximize expected benefits) in decisions made today and tomorrow in the face of uncertainties faced tomorrow. MSPs are an extension of 2-SPs, which involves making decisions sequentially based on new information that emerges over time. For an overview of \ac{SP} and \ac{MSP} models and their applications, we refer the reader to \citet{pereira_1991,goh_2007,zhang_2017, li_2021}. 

MSP models pose significant computational challenges that stem from the inherent uncertainty in the data and the intricate nested structure characteristic of multistage decision-making problems. A common approach to handling uncertainty is to approximate the underlying stochastic process using a scenario tree. A scenario tree branches at each stage to consider a finite number of possible scenarios. When the number of scenarios is small, the deterministic equivalent of the MSP problem can be solved using commercial optimization solvers, such as CPLEX or Gurobi. Lagrangian and Benders decomposition methods are often used to solve problems with a large number of scenarios. These methods break down the problem into smaller subproblems that are easier to solve \citep{laporte_1993, guignard_2003}. Rolling-horizon and heuristic approaches are also used to solve larger-scale MSP models \citep{christian_2015}.

\ac{MSP} models are often used for resource allocation in the healthcare sector, as the size and timing of the demand for these resources are often uncertain, and the availability of resources is limited. For example, \cite{yin_2021b} develop an \ac{MSP} model to effectively manage resource allocation during the 2018-2020 Ebola outbreak. In response to the COVID-19 pandemic, \cite{yin_2023} propose an \ac{MSP} model for allocating ventilators across counties in New York and New Jersey. Their MSP incorporates a compartmental model to account for disease dynamics under varying levels of uncertainty. \cite{hosseini_2023} introduce a multi-stage fuzzy \ac{SP} approach for patient allocation across healthcare facilities during COVID-19.  \cite{tanner_2008} and \cite{yarmand_2014}  propose \ac{SP} models to determine optimal vaccination strategies. \cite{mehrotra_2020} propose an \ac{SP} model for ventilator allocation during COVID-19, focusing on cost minimization. These publications highlight the relevance of MSPs to support sequential decisions related to resource allocation under uncertainty.

Integrating compartmental models within \ac{MSP} resource allocation models provides a robust framework for planning responses to disease outbreaks. This approach enables high-quality decision making over time, taking into account the unpredictable nature of disease spread and resource availability. While some researchers are investigating this modeling approach, there is more to be done. Our proposed data-driven compartmental MSP model addresses \ac{VH} and fairness in resource allocation within this modeling framework. Our proposed model supports effective vaccination strategies and resource allocations that ensure greater immunity and mitigate the impact of outbreaks.

\subsection{Equity in Healthcare Resource Allocation}

A \ac{WHO} report highlights the importance of following ethical principles in the allocation of healthcare resources by analyzing the utility, efficiency, and fairness approaches taken by decision makers \citep{who_equity}. Utilitarian decision makers implement strategies that yield the highest overall benefit or the lowest overall cost. Efficient decision makers advocate using the least amount of resources to achieve the maximum possible impact. Fair decision makers aim to ensure equitable allocation and reduce disparities between groups and individuals. Aligning resource allocation strategies with these principles can be challenging as each strategy, while meeting one principle, falls short of meeting others. The work of \cite{lane_2017} highlights the complexity of ensuring equity in allocating healthcare resources, emphasizing the need for clear definitions that reflect societal values to ensure equal access to healthcare, address patient needs effectively, and distribute potential benefits fairly. Equitable strategies can build public trust and encourage cooperation, which are essential during health crises. 

Several OR models for resource allocation take a utilitarian approach to optimize logistics and operations to control disease spread \citep{sun_2014,liu_2015,liu_2016,zaric_2001,day_2020}. Other studies introduce fairness metrics and propose models that optimize fairness in resource allocation. For example, \cite{orgut_2016} proposes a model that allocates and distributes food donations proportionally to demand, aiming for an equitable distribution while minimizing waste. Similarly, \cite{davis_2015} propose a multi-step optimization approach to improve fairness in the allocation of organs (kidneys) across the U.S., thereby improving equity of organ transplant. \cite{donmez_2022} presents a multi-objective, multi-period, nonlinear model for equitable \ac{PPE} allocation to health centers in pandemics, with the aim of minimizing infections and ensuring fair distribution amidst resource scarcity. The research by \cite{yin_2021} introduces an MSP model for a fair allocation of resources, such as beds and treatment centers, during the West African Ebola outbreak, balancing efficiency and equity in epidemic responses.

Similarly to the literature, our proposed model allows decision makers to evaluate the trade-offs between fairness and efficiency in allocating critical healthcare resources. Unique to our approach is the incorporation of uncertainty in \ac{VH}, which introduces additional complexity to resource allocation strategies. We analyze the `price of fairness', the efficiency loss incurred when striving for equitable resource allocation over time and across different geographic regions. The MSP compartmental model adapts to the evolving landscape of healthcare needs, ensuring that resource allocation remains responsive and equitable even under unpredictable changes in public sentiment towards vaccination.

\section{Method} \label{s:method}

We propose a compartmental \ac{MSP} model to evaluate how uncertainty in \ac{VH}, and disease progression over time and space impact decisions about allocating healthcare resources. Section~\ref{s:method.uncertainty} focuses on the generation of scenarios used in the compartmental MSP model. Detailed descriptions of the compartmental model and the compartmental \ac{MSP} model are presented in Sections~\ref{s:method.compartmental} and \ref{s:method.compartmental_msp}, respectively.

\subsection{\emph{Uncertainty representation}} \label{s:method.uncertainty}

\noindent {\bf Motivation:} Our proposed model considers unccertainty in \ac{VH}. This is motivated by observations made from the \ac{CTIS} conducted by the Delphi Group at Carnegie Mellon University, in collaboration with Facebook \citep{salomone_2021}, and from data about vaccination rates collected by the Centers for Disease Control and Prevention \cite{cdc_2021_county}. 

\begin{figure}[htp]
    \centering
    \includegraphics[width=\textwidth]{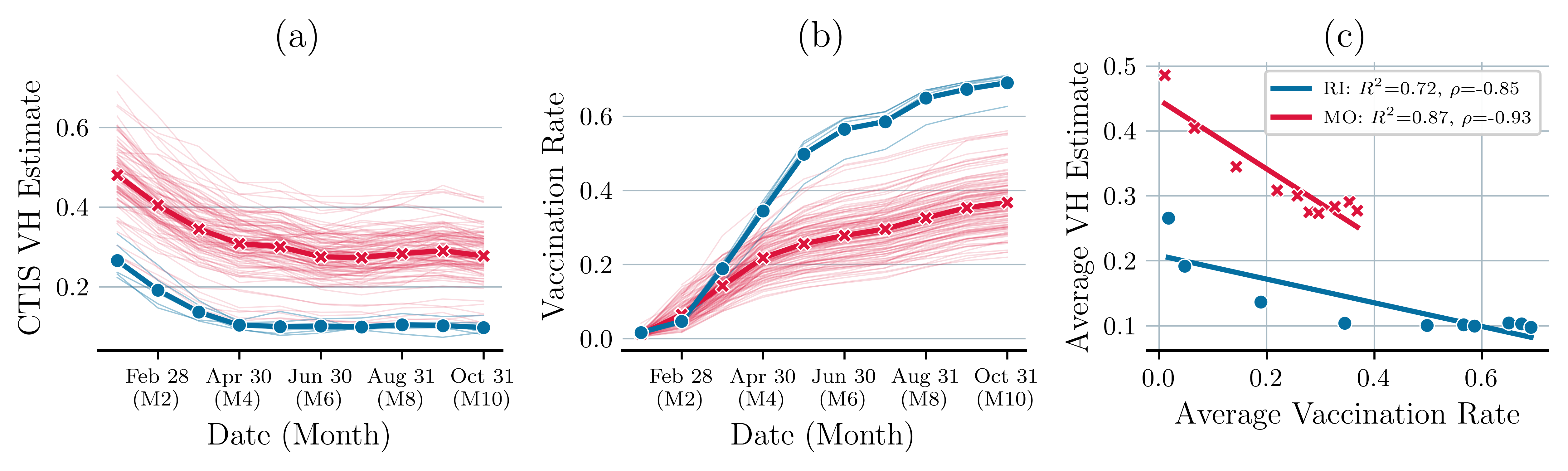}
    \caption{(a) \ac{CTIS}-derived vaccine hesitancy, (b) vaccination uptake trends, and (c) correlation plot between \ac{VH} and vaccination rate of counties in Rhode Island (RI) and Missouri (MO).\label{fig:vh}}
\end{figure}

Figure~\ref{fig:vh}(a) presents \ac{VH} estimates of every county in the states of Missouri (MO) and Rhode Island (RI), providing a snapshot of the willingness of the public to vaccinate. The VH estimate is derived from responses to the following question in the CTIS survey: \textit{``If a vaccine to prevent COVID-19 was offered to you today, would you choose to get vaccinated?''}. Participants were asked to select one of the following answers \textit{``Yes, definitely,''} \textit{ ``Yes, probably,''} \textit{``No, probably not,''} or \textit{``No, definitely not.''} The last three responses are considered indicative of VH. Each thin line represents \ac{VH} progression of an individual county. The thicker line depicts the average \ac{VH} trend for the respective state. 

Figure~\ref{fig:vh}(b) shows the vaccination rates in counties of Missouri and Rhode Island based on data from \cite{cdc_2021_county}. The data show fully vaccinated individuals who have received the complete dose - either two doses of a two-dose vaccine or one dose of a single-dose vaccine. These data highlight a general decline in \ac{VH} over time and underscore the varied vaccination sentiments in counties. Figure~\ref{fig:vh}(c) illustrates a negative correlation between average vaccination rates and average \ac{VH} estimates with correlation coefficients of $\rho = -0.85$ for Rhode Island and $\rho=-0.93$ for Missouri. The high correlation coefficients in both states indicate that higher vaccination rates are associated with lower \ac{VH}. 

Within a state, varying attitudes towards vaccination and differing vaccination rates significantly impact disease transmission and the resulting need for healthcare resources.
We capture variations in VH over time and space via a scenario tree that we describe next.   

\noindent {\bf Scenario tree:} Scenario trees provide a systematic approach to represent uncertainty that is suitable for quantitative models. They provide likely scenarios of the future with associated probabilities. 

Each scenario $\omega$ ($\omega\in\Omega$) represents the changes of VH for a specific population group, over time. We derive temporal changes in \ac{VH} from the historical data provided by the \ac{CTIS}. We fit a distribution to the rate of change of VH. Numerical analysis shows that the (monthly) rate of change of VH is normally distributed. The corresponding mean ($\mu_t$) and standard deviation ($\sigma_t$) change from one period (month) to the next (see Figure~\ref{fig:modeling_vh}). To construct a decision tree of manageable size, we approximate the normal distribution with a discrete triangular distribution that assumes the values $\mu_t, \mu_t - \sigma_t, \mu_t + \sigma_t$ with probabilities of 0.158, 0.684, and 0.158, respectively.

\begin{figure}[htp]
    \centering
    \includegraphics[width=\textwidth]{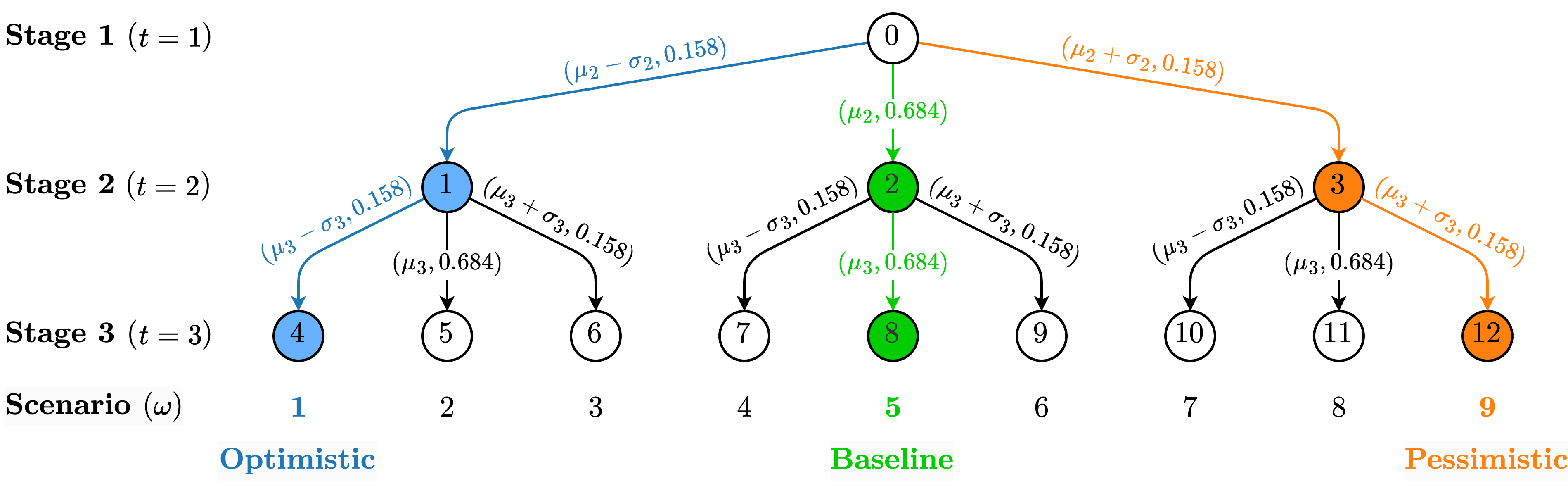}
    \caption{A scenario tree illustrating three stages, with each stage branching into three paths, resulting in nine unique scenarios. Branch annotations denote the rate of change and corresponding probabilities.\label{fig:tree}}
\end{figure}

Figure~\ref{fig:tree} presents the scenario tree for a problem with $T=3$ stages. The resulting scenario tree has $3^2=9$ unique scenarios. Let $p_{\omega}$ denote the probability associated with each of these scenarios. The tree showcases a baseline (green path), an optimistic (blue path), and a pessimistic (orange path) scenario. In the optimistic scenario, the VH changes by $\mu_t - \sigma_t$ in each stage. In the pessimistic scenario, the VH changes by $\mu_t + \sigma_t$ at each stage. In the baseline scenario, the VH changes by $\mu_t$ at each stage. Consequently, the probability associated with the pessimistic and optimistic scenarios is $0.158^2=0.024$. The probability associated with the baseline scenarios is $0.684^2=0.468$.

\subsection{\emph{A compartmental model of  COVID-19}} \label{s:method.compartmental}

The \ac{SEIR} model divides the population into four groups (compartments): Susceptible (\(\pmb{S}\)), Exposed (\(\pmb{E}\)), Infectious (\(\pmb{I}\)), and Recovered (\(\pmb{R}\)). The model captures the transition of individuals through these compartments over time. It also accounts for the latency period before an individual becomes infectious, which is important in modeling the spread of diseases such as COVID-19. In a previous study, our team extended this model to address the specific complexities of COVID-19, such as vaccine availability, VH, and the scarcity of critical healthcare resources in the spread of the disease \citep{BER_24}. For the convenience of the reader, we summarize the model in the following paragraphs.  

\begin{figure}[htp]
    \centering
    \includegraphics[width=\textwidth]{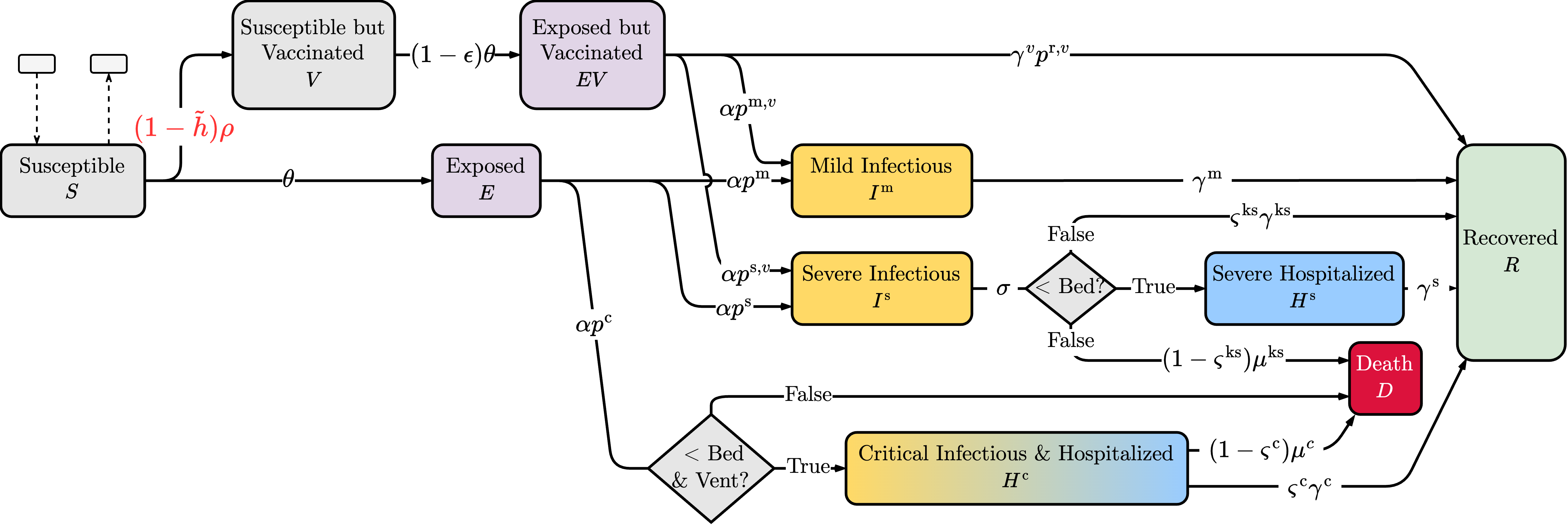}
    \caption{\ac{SVEIHR}: COVID-19 compartmental model that incorporates vaccination dynamics, VH, and healthcare resource constraints.\label{fig:sveihr}}
\end{figure}

A schematic representation of the proposed \ac{SVEIHR} model is shown in Figure~\ref{fig:sveihr}. In this model, susceptible individuals can be exposed or opt to vaccinate, moving forward to the vaccinated compartment (\pmb{V}). The size of compartment \pmb{V} is dictated by $\rho$, representing the maximum vaccination rate achievable in the absence of \ac{VH}. The willingness of people to get vaccinated is represented by $(1-\tilde{h})$, where, $\tilde{h}$ is a stochastic parameter that represents the proportion of the population reluctant to vaccinate. Furthermore, the model considers that vaccination does not confer absolute immunity; therefore, vaccinated individuals can progress to an Exposed but Vaccinated compartment ($\pmb{EV}$), with a probability determined by vaccine efficacy, $\epsilon$.

The model captures different postexposure health trajectories based on incubation rates $\alpha$ and the severity of the disease. Individuals progress to Mild ($\pmb{I^{\text{m}}}$), Severe ($\pmb{I^{\text{s}}}$), or Critical Condition ($\pmb{H^{\text{c}}}$) compartments with probabilities $p^{\text{m}}, p^{\text{s}},$ and $p^{\text{c}}$, respectively. Vaccinated individuals are less likely to experience severe/critical health outcomes; thus, several infected individuals are likely to recover and move to the compartment ($\pmb{R}$). 

The model factors the availability of critical healthcare resources, such as hospital beds (\textbf{Bed}) and ventilators (\textbf{Vent}), in determining patients' health outcomes. Refer to the example presented in Appendix~\ref{appendix:compartmental_model} for a demonstration of the impact of the number of beds and ventilators on the model outcomes. Recovery ($\gamma$) and mortality ($\mu$) rates vary depending on the severity of the disease and the availability of healthcare resources. The model accounts for the transition of critically ill patients to the compartment ($\pmb{D}$) when ventilators are unavailable, emphasizing the dire consequences of resource limitations during a pandemic.

\subsection{\emph{A compartmental multi-stage stochastic optimization model}} \label{s:method.compartmental_msp}

This section outlines the formulation of the proposed  \acf{MSP} model. The model determines an optimal sequence of decisions related to ventilator allocation across various regions over time. The \ac{MSP} minimizes the total number of deaths by dynamically adjusting the distribution of healthcare resources in response to changing VH rates and state transitions within the \ac{SVEIHR} model. Resource allocation decisions take into account not just immediate health outcomes, but also the need for equitable access to care. By integrating the compartmental model with the MSP model, we now better understand how \ac{VH} and disease dynamics affects healthcare systems, and use this information to optimize the allocation of ventilators. We refer to this proposed model as the compartmental MSP model. For a detailed description of the model notations used, please see Appendix~\ref{appendix:notation}.

\noindent{\bf Model assumptions:} \label{s:method.model.assumptions}
Our model operates on three primary sets of assumptions that are critical to its structure and function:

\noindent \emph{Assumption 1: Uncertainty in vaccine hesitancy (VH):} The model considers three potential outcomes of \ac{VH} in each decision stage. Let ${h}_{t-1}$ represent the proportion of the population reluctant to vaccinate in period $t-1$. Then, in period $t$, the proportion of the population reluctant to vaccinate for the optimistic, expected, and pessimistic scenarios are $h_{t-1}(1 + \mu_t-\sigma_t); h_{t-1}(1 + \mu_t)$; and $h_{t-1}(1 + \mu_t + \sigma_t)$. Each of these outcomes is associated with a fixed probability that does not change over the course of the study. These probabilities remain constant and unaffected by changing decisions. We note that decision makers can customize these outcomes (number of outcomes and values) based on available data. 

\noindent \emph{Assumption 2: Disease dynamics}: Since \ac{SVEIHR} is a compartmental epidemiological model, the following assumptions of \ac{SEIR} model apply. ($i$) Every individual in a population is equally likely to interact with others. As a result, disease progression and transitions between compartments are consistent throughout the population. However, it should be noted that in a real-world setting, elements such as social networks and age demographics can significantly influence disease transmission patterns. ($ii$) The size of the total population does not change due to natural births, non-COVID-19-related deaths, or migration. Since the time frame considered in our study spans several months, recovered individuals are assumed to have full and lasting immunity to the disease. ($iii$) Infection, recovery, and exposure rates are assumed to be deterministic. However, we recognize that in practice, they can fluctuate for many reasons, from policy interventions and demographic differences to viral mutations.

\noindent \emph{Assumption 3: Healthcare system:} The available healthcare resources (i.e., beds, ventilators) are only used to treat COVID-19 patients. Furthermore, our proposed MSP model does not account for potential reallocations of healthcare resources.

\noindent{\bf Mathematical model formulation:}\label{s:method.model.formulation}
Classical SIR are continuous-time models. However, several researchers propose discrete-time approximations of these models that are easier to solve \citep{wacker2020time}. These studies show that many of the desired properties of the time-continuous SIR models remain valid in the time-discrete approximations. We use a discrete-time model to mimic the dynamic changes in the size of compartments of the \ac{SVEIHR} model over time.  Let $\mathcal{T}$ represent the set of time periods in the planning horizon. Let $\mathcal{J}$ represent the set of decision periods (stages). Note that, between every two consecutive decision periods, there are several time periods, thus $\mathcal{J} \subset \mathcal{T}$. We denote by $\mathcal{R}$ the set of regions studied. We distinguish between a decision stage $j$ and a time period $t$ to keep the discrete approximation of the SEIR model effective, ensuring short time intervals between transitions, while avoiding too many decision stages for the MSP model. 

In a decision period (stage), we determine the allocation of ventilators. Our decision variables are $x_{\omega,j,r}$, which represent the number of ventilators assigned to the region $r\in \mathcal{R}$ at the beginning of period $j\in\mathcal{J}$ in scenario $\omega\in \Omega$. Let $\mathcal{X}_{\omega,t,r}$ represent the cumulative number of ventilators allocated to region $r$ by time $t$.

The decisions made affect the state variables, which are the number of susceptible individuals who are unvaccinated ($S_{\omega,t,r}$) and vaccinated ($V_{\omega,t,r}$); the number of exposed individuals who are unvaccinated ($E_{\omega,t,r}$) and vaccinated ($EV_{\omega,t,r}$); the number of infectious  individuals with mild ($I^{\text{m}}_{\omega,t,r}$) and and severe ($I^{ \text{s}}_{\omega,t,r}$) symptoms; the number of hospitalized individuals in severe ($H^{\text{s}}_{\omega,t,r}$) and in critical ($H^{\text{c}}_{\omega,t,r}$) condition; the number of individuals admitted to a hospital and in severe ($\mathcal{A}^{\text{s}}_{\omega,t,r}$) or critical ($\mathcal{A}^{\text{c}}_{\omega,t,r}$ ) condition; the number of individuals in severe ( $\mathcal{K}^{\text{s}}_{\omega,t,r}$) and critical ($\mathcal{K}^{\text{c}}_{\omega,t,r}$) condition not admitted to the hospital due to capacity constraints; the number of recovered individuals ($R_{\omega,t,r}$); and the number of deaths ($D_{\omega,t,r}$). Let $\mathcal{L} = \{S, V, E, EV, I^{\text{m}}, I^{\text{s}}, H^{\text{s}}, H^{\text{c}}, R, D\}$  represent the set of our state variables.

Next, we provide detailed descriptions of the objective function and constraints of the proposed compartmental MSP.

\noindent \emph{Objective function:} The objective is to minimize the total expected number of deaths across regions during the planning horizon. 
\begin{equation}
  \mathcal{Z} = \min \sum_{r\in \mathcal{R}}\sum_{\omega \in \Omega} p_{\omega}   D_{\omega,|\mathcal{T}|,r} \label{eq:obj}
\end{equation}

\noindent \emph{Constraints:} We organize our constraints into five groups {(\textit{G1, G2,\ldots, G5})}, each tailored to fulfill a specific function. This structure enhances the model's clarity and highlights the essential contribution of each constraint category to the model's objective.

\noindent (\textit{G1:}) \emph{\ac{SVEIHR} compartmental model:} This group of constraints captures the dynamic changes of state variables over time. We begin by initializing the state variables $\ell\in\mathcal{L}$ at $t=0$. Let $\boldsymbol{\pi}_{\ell,r}$ represent the corresponding initial values. 

\begin{equation}
\ell_{\omega,0,r} = \boldsymbol{\pi}_{\ell,r} \hspace{0.1in} \forall \ell\in\mathcal{L}, r\in\mathcal{R}, \omega\in\Omega\label{eq:init}
\end{equation}

Constraints~\eqref{eq:susceptible}-\eqref{eq:death} describe the transition of the population between different compartments of the epidemiological model over time (as illustrated in Figure ~\ref{fig:sveihr}). Constraint~\eqref{eq:susceptible} updates the size of the susceptible population in compartment \pmb{S}. To determine the size of \pmb{S} at the end of period $t+1$, we update the size of \pmb{S} at the end of period $t$ by subtracting the number of individuals exposed, vaccinated and subtracting (adding) those who immigrated (emigrated) to the region from (to) other regions considered in the analysis. In Appendix~\ref{appendix:migration_flow}, we provide descriptions of the effects of migration in our model. Constraint~\eqref{eq:vaccinated} updates the size of the vaccinated population, compartment \pmb{V}. To determine the size of \pmb{V} at the end of the period $t+1$, we update the size of \pmb{V} at the end of period $t$ by subtracting the number of individuals exposed and adding the number of individuals vaccinated. The efficacy rate of the vaccine, $\epsilon$, affects the size of this compartment. Constraint~\eqref{eq:exposed} updates the size of the unvaccinated population who were exposed, compartment \pmb{E}. To determine the size of  \pmb{E} at the end of period $t+1$, we update the size of \pmb{E} at the end of period $t$ subtracting the number of individuals exposed who became infected, and by adding the number of individuals exposed during period $t+1$. Constraint~\eqref{eq:ev} updates the size of the vaccinated population who are exposed, compartment \pmb{EV}. This constraint updates the size of \pmb{EV} by taking a similar approach to ~\eqref{eq:exposed}. The difference between constraints~\eqref{eq:exposed} and~\eqref{eq:ev} is that vaccinated individuals can recover after exposure and do not become critically ill. Constraints~\eqref{eq:Im} and~\eqref{eq:Is} update the size of compartments \pmb{$I^{\text{m}}$} and \pmb{$I^{\text{c}}$}, respectively. To determine the size of compartment \pmb{$I^{\text{m}}$} at the end of $t+1$, constraint~\eqref{eq:Im} considers new infections, and recoveries. To determine the size of compartment \pmb{$I^{\text{s}}$} at the end of $t+1$, constraint~\eqref{eq:Is} considers new infections, recoveries, deaths and transitions to the hospitalization compartment. Constraints~\eqref{eq:h_s} and~\eqref{eq:h_c} update the size of compartments  \pmb{$H^{\text{s}}$} and \pmb{$H^{\text{c}}$}. To determine the number of severely hospitalized individuals (the size of compartment \pmb{$H^{\text{s}}$}) at the end of period $t+1$, constraint~\eqref{eq:h_s} considers new admissions, recoveries, and deaths. To determine the number of critically hospitalized individuals (the size of compartment \pmb{$H^{\text{c}}$}) at the end of period $t+1$, constraint~\eqref{eq:h_c} also considers new admissions, recoveries, and deaths. Constraint~\eqref{eq:r} determines the number of recovered individuals (the size of the compartment \pmb{R}) at the end of the period $t+1$ by aggregating the number of those recovered from various compartments (of infected and hospitalized individuals), each with specific recovery rates. Lastly, the constraint~\eqref{eq:death} determines the size of the compartment \pmb{D}. It compiles the death toll within the region, including deaths from the inability to hospitalize individuals in severe and critical conditions, and deaths among individuals in critical conditions, adjusting to their respective mortality rates.

\noindent Constraints~\eqref{eq:susceptible}-\eqref{eq:death} are for all $t \in \{1,\ldots,|\mathcal{T}|-1\}, r\in \mathcal{R}, \omega \in \Omega$:

\begin{align}
    S_{\omega, t+1,r} &= S_{\omega,t,r} - \theta_{\omega,t,r}S_{\omega,t,r} -\rho_{r}(1-h_{\omega,t,r})S_{\omega,t,r} + \sum_{r'\in \mathcal{R}} \nu_{r'\rightarrow r} S_{\omega,t,r'} - \sum_{r'\in \mathcal{R}} \nu_{r\rightarrow r'} S_{\omega,t,r}\label{eq:susceptible}\\
    V_{\omega,t+1,r} &= V_{\omega,t,r} + \rho_{r}(1-h_{\omega,t,r})S_{\omega,t,r} - (1-\epsilon)\theta_{\omega,t,r}V_{\omega,t,r} \label{eq:vaccinated} \\
    E_{\omega,t+1,r} &= E_{\omega,t,r} + \theta_{\omega,t,r}S_{\omega,t,r} - \alpha_r E_{\omega,t,r} \label{eq:exposed}\\
    EV_{\omega,t+1,r} &= EV_{\omega,t,r} + (1-\epsilon)\theta_{\omega,t,r}V_{\omega,t,r} - \left(\alpha_r (p^{\text{m},v}_r+p^{\text{s},v}_r) + \gamma^v_r p^{\text{r},v}_r\right)  EV_{\omega,t,r} \label{eq:ev} \\
    I^{\text{m}}_{\omega,t+1,r} &= I^{\text{m}}_{\omega,t,r} + \alpha_r p^{\text{m}}_{r} E_{\omega,t,r} + \alpha_r p^{\text{m},v}_{r} EV_{\omega,t,r} - \gamma^{\text{m}}_r I^{\text{m}}_{\omega,t,r} \label{eq:Im} \\
    I^{\text{s}}_{\omega,t+1,r} &= I^{\text{s}}_{\omega,t,r} + \alpha_r p^{\text{s}}_{r} E_{\omega,t,r} + \alpha_r p^{\text{s},v}_{r} EV_{\omega,t,r} - \sigma_r I^{\text{s}}_{\omega,t,r} \label{eq:Is} \\
    H^{\text{s}}_{\omega,t+1,r} &= H^{\text{s}}_{\omega,t,r} + \mathcal{A}^{\text{s}}_{\omega,t,r} - \gamma^{\text{s}}_r H^{\text{s}}_{\omega,t,r} \label{eq:h_s} \\
    H^{\text{c}}_{\omega,t+1,r} &= H^{\text{c}}_{\omega,t,r} + \mathcal{A}^{\text{c}}_{\omega,t,r} - (1-\varsigma^{\text{c}}_r) \mu^{\text{c}}_r H^{\text{c}}_{\omega,t,r} - \varsigma^{\text{c}}_r \gamma^{\text{c}}_r H^{\text{c}}_{\omega,t,r} \label{eq:h_c} \\
    R_{\omega,t+1,r} &= R_{\omega,t,r} +\gamma^v_r p^{\text{r},v}_{r} EV_{\omega,t,r} + \gamma^{\text{m}}_r I^{\text{m}}_{\omega,t,r} + \gamma^{\text{s}}_r H^{\text{s}}_{\omega,t,r} + \varsigma^{\text{ks}}_r \gamma^{\text{ks}}_r \mathcal{K}^{\text{s}}_{\omega,t,r} + \varsigma^{c}_r\gamma^{\text{c}}_r H^{\text{c}}_{\omega,t,r} \label{eq:r} \\
    D_{\omega,t+1,r} &= D_{\omega,t,r} + (1-\varsigma^{\text{ks}}_r) \mu^{\text{ks}}_r \mathcal{K}^{\text{s}}_{\omega,t,r} + (1-\varsigma^{\text{c}}_r)\mu^{\text{c}}_r H^{\text{c}}_{\omega,t,r} + \mathcal{K}^{\text{c}}_{\omega,t,r} \label{eq:death}
\end{align}

\noindent(\textit{G2}) \emph{Resource capacity constraints:} Constraints~\eqref{eq:cap_vent} to~\eqref{eq:cumm_x_pass} determine the values of the decision variables. The initial number of ventilators in hospitals is limited. Consider that $\Delta_{j}$ is the number of new ventilators available for allocation in period (stage) $j \in \mathcal{J}$. Constraints~\eqref{eq:cap_vent} determine the allocation of these ventilators to hospitals in different regions. Constraints \eqref{eq:cumm_x} determine the total number of ventilators available in region $r$ at period $j$. Since the number of time periods between every two consecutive decision periods is greater than one, we ensure (via constraints~\eqref{eq:cumm_x_pass}) that the total number of available ventilators does not change until the next decision period.    

\begin{align}
    &\sum_{r\in \mathcal{R}} x_{\omega,j,r} \leq \Delta_{j}  && \forall j\in \mathcal{J}, \omega \in \Omega, \label{eq:cap_vent}\\
    &\mathcal{X}_{\omega,j,r} = \mathcal{X}_{\omega,j-1,r} + x_{\omega,j,r} && \forall j \in \mathcal{J}, r\in\mathcal{R}, \omega \in \Omega, \label{eq:cumm_x} \\
    &\mathcal{X}_{\omega,t,r} = \mathcal{X}_{\omega,t-1,r} && \forall t \in \mathcal{T} \setminus \mathcal{J}, r\in\mathcal{R}, \omega \in \Omega. \label{eq:cumm_x_pass} 
\end{align}

Constraint~\eqref{eq:admit_critical} connects the decision variables with the state variables. This constraint determines the number of critically ill patients admitted to a hospital, $\mathcal{A}^{\text{c}}_{\omega,t,r}$. This number depends on the demand for critical care, the number of ventilators available, and the number of beds available. 

Constraint~\eqref{eq:admit_severe} determines the number of severely ill patients admitted to a hospital, $\mathcal{A}^{\text{s}}_{\omega,t,r}$. This number depends on the demand for care and beds available after the admission of patients in critical condition. In Appendix~\ref{appendix:linearize} we detail the linearization of ``min" functions in constraints~\eqref{eq:admit_critical} and~\eqref{eq:admit_severe}. The restrictions~\eqref{eq:deny_critical} and \eqref{eq:deny_severe} determine the number of patients with severe and critical conditions who were not admitted to a hospital due to limited resources. Constraints~\eqref{eq:admit_critical} to~\eqref{eq:deny_severe} are for all $t \in \mathcal{T}, r\in \mathcal{R}, \omega \in \Omega$.

\begin{align}
    \mathcal{A}^{\text{c}}_{\omega,t,r} &= \min\Bigl\{\alpha_r p^{\text{c}}_{r} E_{\omega,t,r},\: \mathcal{X}_{\omega,t,r} - H^{\text{c}}_{\omega,t,r}, b_{r} - H^{\text{c}}_{\omega,t,r} - H^{\text{s}}_{\omega,t,r} \Bigr\} \label{eq:admit_critical}\\
    \mathcal{A}^{\text{s}}_{\omega,t,r} &= \min\Biggl\{\sigma I^{\text{s}}_{\omega,t,r},\: b_{r} - H^{\text{c}}_{\omega,t,r} - H^{\text{s}}_{\omega,t,r} - \mathcal{A}^{\text{c}}_{\omega,t,r} \Biggr\} \label{eq:admit_severe}\\
    \mathcal{K}^{\text{c}}_{\omega,t,r} &= \alpha_r p^{\text{c}}_{r} E_{\omega,t,r} - \mathcal{A}^{\text{c}}_{\omega,t,r} \label{eq:deny_critical}\\
    \mathcal{K}^{\text{s}}_{\omega,t,r} &= \sigma I^{\text{s}}_{\omega,t,r} - \mathcal{A}^{\text{s}}_{\omega,t,r} \label{eq:deny_severe}
\end{align}

\noindent (\emph{G3}:)  \emph{Resource allocation constraints:} Equity in healthcare resource distribution, especially during emergencies, is vital for public health policy. Defining and measuring equity in public health decisions is a complex task with no widely accepted standard. We present several approaches used by the literature to ensure fair allocation of healthcare resources (see our discussion of \eqref{eq:equity}-\eqref{eq:equal_alloc}). Each constraint, when applied individually, offers a distinct approach to resource allocation. Constraint~\eqref{eq:equity}, adapted from \citep{yin_2021}, uses information about the proportion of the expected critically ill individuals in a region (the first term in this constraint) and the corresponding population proportion of the region (the second term in this constraint) to allocate ventilators. It uses an equity threshold parameter, denoted by $k \in (0, 1]$, to gauge need-based equity in resource allocation. As $k$ approaches zero, equity in resource allocation is enforced. As the value of $k$ approaches one, the equity requirements become less stringent.  In Appendix~\ref{appendix:linearize}, you will find our approach to linearizing this constraint.  

\begin{align}
    \left\lvert\frac{\sum_{t \in \mathcal{T}}\sum_{\omega\in\Omega}p_\omega H^{\text{c}}_{\omega,t,r}}{\sum_{t \in \mathcal{T}}\sum_{r'\in \mathcal{R}}\sum_{\omega\in\Omega} p_\omega H^{\text{c}}_{\omega,t,r'} } - \frac{n_r}{\sum_{r'\in \mathcal{R}}n_{r'}} \right\rvert \leq k \label{eq:equity}
\end{align}

The restrictions~\eqref{eq:proportional} present a variation in the population-proportional method for the allocation of resources \citep{dangerfield_2019, who_2020_equity}. We include an equity threshold parameter $\zeta$. This constraint ensures that each region receives ventilators in proportion to its share of the total population, with $\zeta$ providing adjustable scaling to fine-tune how strictly the allocation of resources follows population size. For $\zeta = 1$, the additional ventilators available for use during period $j$, denoted as $\Delta_j$, are distributed according to the population size of each region. As $\zeta$ decreases, this constraint relaxes, and the equity requirements become less stringent.

\begin{align}
    \sum_{\omega \in \Omega} p_\omega x_{\omega,j,r} \geq \left(\frac{n_r}{\sum_{r' \in \mathcal{R}} N_{r'}}\right) \zeta \Delta_{j}, \hspace{0.1in} \forall j\in\mathcal{J}, r\in\mathcal{R}.\label{eq:proportional}
\end{align}

Both constraints \eqref{eq:equity} and \eqref{eq:proportional}, focus on ensuring an equitable distribution of critical resources. One of the constraints uses a need-based approach, and the other population-based approach to resource allocation.

Constraints~\eqref{eq:equal_alloc} ensure that ventilators are distributed equally, ignoring differences in population size or the specific needs of each region. As a result, the same number of ventilators is allocated to every region. 

\begin{align}
    \sum_{\omega \in \Omega} p_\omega x_{\omega,j,r} \geq \left\lfloor \dfrac{\Delta_j}{|\mathcal{R}|}\right\rfloor, \hspace{0.1in} \forall j \in \mathcal{J}, r \in \mathcal{R}. \label{eq:equal_alloc}
\end{align}

Notice that we solve the compartmental MSP model with either constraints~\eqref{eq:equity}, or~\eqref{eq:proportional}, or \eqref{eq:equal_alloc} to enforce equity or equality in resource allocation. We also solve the compartmental MSP model without constraints ~\eqref{eq:equity},  ~\eqref{eq:proportional}, and~\eqref{eq:equal_alloc} to demonstrate the effects of a utilitarian approach that aims to minimize the overall cost of resource allocation.

\noindent(\textit{G4}) \emph{Non-anticipativity constraints (NACs):} These constraints ensure that decisions made at a given time (decision stage) are based solely on the information available up to that point. This prevents the use of future, yet unknown, information in decision making. Specifically, constraints \eqref{eq:NACs} ensure that decisions that share an identical scenario path up to time $j$ must be consistent and identical across shared paths. In \eqref{eq:NACs}, $\mathcal{S}_{\omega,j,r}$ represents the set of scenarios that are indistinguishable from $\omega$ at time $t$. In Appendix~\ref{appendix:nac}, we provide an example that demonstrates the implementation of these constraints.

\begin{align}
    x_{\omega,j,r} - x_{\omega',j,r} &= 0, &\forall \omega \in \Omega, \omega' \in \mathcal{S}_{\omega,j,r}, j\in\mathcal{J}, r\in\mathcal{R}. \label{eq:NACs}
\end{align}

\noindent(\textit{G5}) \emph{Other constraints:} Constraints~\eqref{eq:define_x} restrict the decision variables to be non-negative integers. Constraints~\eqref{eq:define_l} and \eqref{eq:define_A} restrict the state variables to be nonnegative. The following constraints are also for all $\omega \in \Omega, r\in\mathcal{R}$. 

\begin{align}
    &x_{\omega,j,r} \in \mathbb{Z}_{\geq 0} & \forall j \in \mathcal{J},  \label{eq:define_x}\\
    &\ell_{\omega,t,r} \geq 0 & \forall \ell \in \mathcal{L}, t \in \mathcal{T}, \label{eq:define_l}\\
    &\mathcal{A}^{\text{s}}_{\omega,t,r}, \mathcal{A}^{\text{c}}_{\omega,t,r}, \mathcal{K}^{\text{s}}_{\omega,t,r}, \mathcal{K}^{\text{c}}_{\omega,t,r} \geq 0 & \forall t \in \mathcal{T}. \label{eq:define_A}
\end{align}

\section{Data Collection and Model Verification and Validation}

To validate the compartmental MSP model, we developed a case study using county-level data from Arkansas during January to May 2021. During this period of the COVID-19 pandemic, vaccines were available to the public. Next, we discuss our data collection, data processing, and our approach to model validation and verification. The \ac{SVEIHR} and \ac{MSP} compartmental models presented in Section~\ref{s:method} are modeled using Python 3.10. The optimization model in Section~\ref{s:method.compartmental_msp} is solved using the DOcplex API from CPLEX version 22.1.1. We performed our computational analysis on a system equipped with a 32-core CPU and 192 GB of RAM.

\subsection{\emph{Data collection and processing}} \label{s:discussion.data.processing}

\noindent{\bf Data collection:} We collect data about county-level vaccination uptake from January to May 2021 from the \cite{cdc_2021_county}, and data about \Ac{VH} from the CTIS dataset \citep{salomone_2021}. Data on healthcare resources come from the \citet{covidcaremap} project. Data are organized by facility, county, referral region, and state. Table~\ref{tab:healthcare_cap} summarizes the data for Arkansas.

\begin{table}
    \centering
    \caption{Healthcare capacity by region in Arkansas.\label{tab:healthcare_cap}}
    \begin{tabular}{ccccc}
        \toprule
        \textbf{Region} & \textbf{\# Counties} & \textbf{Population (Ratio)} & \textbf{Licensed Beds} & \textbf{ICU Beds} \\
        \midrule
        R1 & 18 & 394,446 (14.45\%) & 1,405 & 77 \\
        R2 & 14 & 1,551,512 (56.84\%) & 9,357 & 566 \\
        R3 & 13 & 171,946 (6.30\%) & 490 & 18 \\
        R4 & 30 & 611,686 (22.41\%) & 2,032 & 105 \\
        \midrule
        \textbf{Total} & 75 & 2,729,590 & 13,284 & 766 \\
        \bottomrule
    \end{tabular}
\end{table}

\begin{figure}[htp]
    \centering
    \includegraphics[width=1\textwidth]{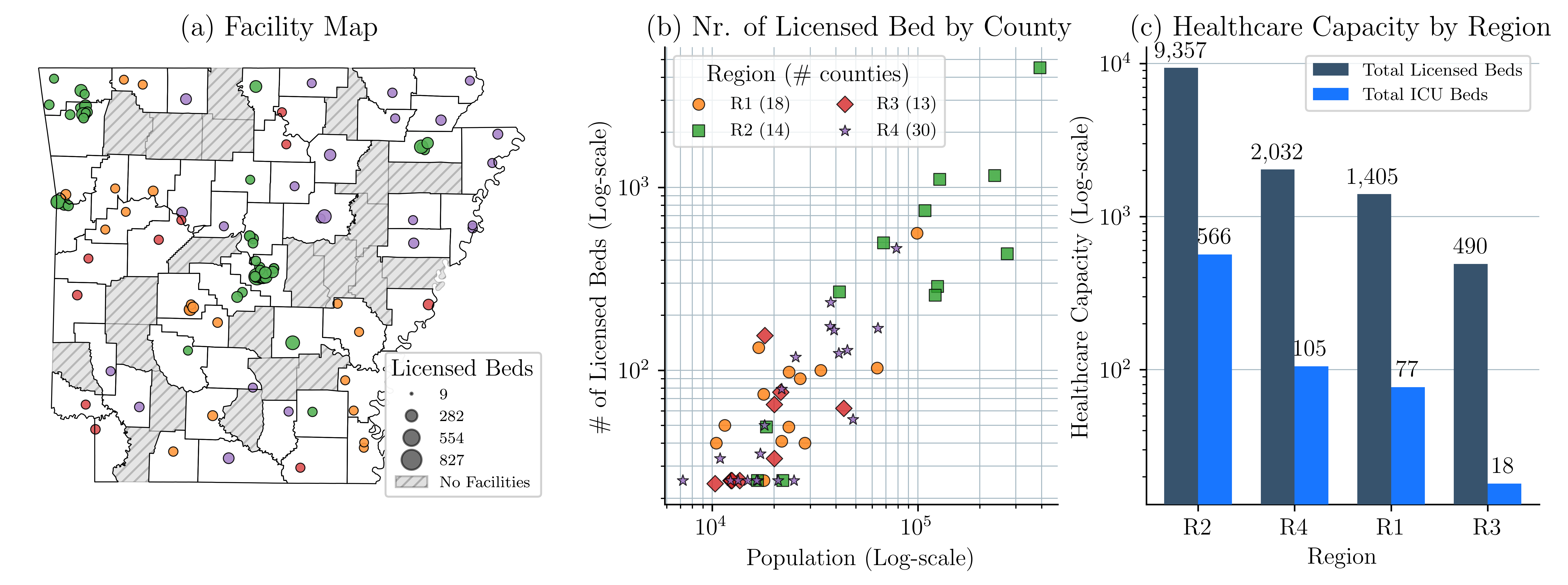}
    \caption{(a) County-level distribution of licensed hospital beds in Arkansas. (b) Relationship between population size and number of licensed hospital beds by county. (c) The total number of licensed and ICU beds, aggregated by region of study.}
    \label{fig:healthcare_cap}
\end{figure}

Figure~\ref{fig:healthcare_cap}(a) presents the county-level distribution of hospitals and licensed beds in Arkansas at the beginning of the COVID-19 pandemic. Figure~\ref{fig:healthcare_cap}(b) illustrates the relationship between population size and the number of licensed beds in each county. We cluster counties of Arkansas into four regions based on their \ac{VH}. Below, we provide details of our clustering approach. Figure~\ref{fig:healthcare_cap}(c) summarizes the total number of licensed beds and \ac{ICU} beds in each region. We assume that each ICU bed has a ventilator. The results of these figures emphasize the variations in healthcare capacity across the state. 

\begin{figure}[htp]
    \centering
    \includegraphics[width=\textwidth]{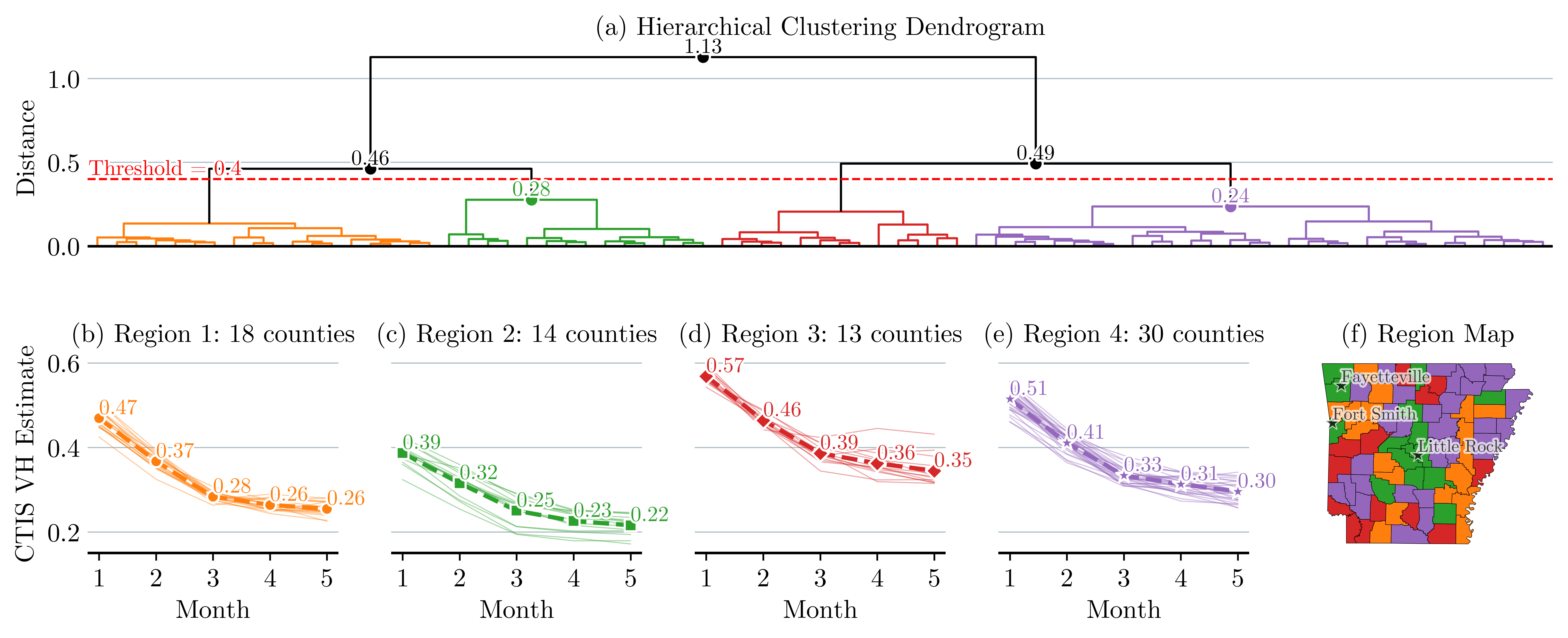}
    \caption{(a) Hierarchical clustering of \ac{CTIS} VH data for counties in Arkansas. (b)-(e) Region-specific VH trends over five months (January to May 2021). (f) Spatial distribution of the four regions.}
    \label{fig:clustering}
\end{figure}

\noindent{\bf Clustering of the data:} In an effort to reduce the size of the decision tree of the \ac{MSP} model, we cluster counties of Arkansas into four regions. The counties are clustered according to their \ac{VH} similarity. We measure similarity using the Euclidean distance between county-level \ac{VH} time series data. We use the SciPy Python dendrogram tool to visualize the clustering of counties that belong to the same region (Figure~\ref{fig:clustering}(a)) \citep{scipy}. We use as a threshold a distance of 0.4 and identify four distinct regions. Figures~\ref{fig:clustering}(b)-(e) illustrate \ac{VH} trends within each region. Each thin line represents the \ac{VH} of a county, and the dotted line represents the average \ac{VH} of that region. In particular, Region 3, marked red, starts with the highest \ac{VH} and, despite the declining trend, remains the highest throughout the study period. In contrast, Region 2, marked green, begins with a \ac{VH} of 0.39, which then drops dramatically, ending as the region with the lowest \ac{VH}. Figure~\ref{fig:clustering}(f) maps the regions geographically, highlighting that Arkansas's major cities, Fayetteville, Fort Smith, and Little Rock, are in Region 2, which has the lowest \ac{VH}.

\noindent{\bf Other data processing:} Our model employs the parameter $\Delta_j$ to indicate the number of ventilators available for distribution among different regions at the start of each decision period, $j$. The COVID-19 pandemic significantly increased the demand for ventilators, a critical medical device that saved many lives. In response, under \acf{DPA}, various manufacturers, including auto manufacturers, quickly began to increase the production of ventilators or repurpose their facilities for this task \citep{albergotti_2020}. For example, \ac{GM} was contracted to supply 30,000 medical ventilators for the national stockpile \citep{cnn_gm_2020}. In April 2020, AdvaMed reported that its member companies were ramping up production from an average of 2,000 to 3,000 ventilators per week to an expected 5,000 to 7,000 per week. This was a notable increase from the 700 ventilators they produced weekly for the U.S. market in 2019 \cite{advemed_2020}. Due to the lack of real-world data related to the additional ventilators allocated to Arkansas after the enactment of \ac{DPA}, we use the following approach to generate these data. 

We use equations~\eqref{eq:vent_calc0} to~\eqref{eq:vent_calc2} to determine $\Delta_j$. Let $\underline{p}$ represent the initial stockpile of ventilators in Arkansas prior to enactment of \ac{DPA}, and let $\overline{p}$ represent an upper bound on the total number of available ventilators. Let $p_j$ represent the additional ventilators made available in period $j$.

\begin{align}
    \Delta_1 &= \underline{p} & \label{eq:vent_calc0} \\
    \Delta_j &= \Delta_{j-1} + p_j & \forall j \in \mathcal{J}\setminus\{1\} \label{eq:vent_calc1} \\
    \Delta_{j} &\leq \overline{p}  & \forall j \in \mathcal{J} \label{eq:vent_calc2}
\end{align}

We begin with an initial stockpile ($\underline{p}$) of 100 ventilators. In each subsequent time period, we increase the number of ventilators ($p_j$) by 50. As a result, the number of ventilators available for distribution in Arkansas follows this pattern: $\Delta_1 = 100$, $\Delta_5 = 150$, $\Delta_9 = 200$, $\Delta_{13} = 250$, and $\Delta_{17} = 300$.

Please refer to tables in Appendix~\ref{appendix:model_data} for a complete list of the parameter values used in our \ac{MSP} compartmental model. The table also presents the tuned parameters, which are outcomes of the model calibration process detailed in Section~\ref{s:model.calibration}.

\subsection{\emph{Calibration of the compartmental model}} \label{s:model.calibration}

The calibration of \ac{SVEIHR} focuses on fine-tuning the model parameters to ensure that it accurately represents the dynamics of disease progression in each region of study. We collected data from COVID-19 Open Data \citep{Wahltinez2020}. This dataset provides county-level data. For the purpose of this study, we aggregate county-level data for Arkansas into the four aforementioned regions. 

We use data from two critical time periods in the U.S. The first spans 20 weeks, from October 25, 2020 to March 7, 2021. This period includes the peak of the wave of coronavirus infections during the Fall of 2020. We use these data to evaluate \ac{SVEIHR}'s power in predicting new infections. The second spans 20 weeks, from January 10 to May 23, 2021, and corresponds to the initial stages of vaccine rollout. We use these data to evaluate \ac{SVEIHR}'s power in predicting vaccinations. These study periods do not coincide with the start of the pandemic. Therefore, we adjust the size of the initial susceptible population to account for individuals who were infected, recovered, or died prior to the start of the study period. 

\begin{figure}[htp]
    \centering
    \includegraphics[width=1\textwidth]{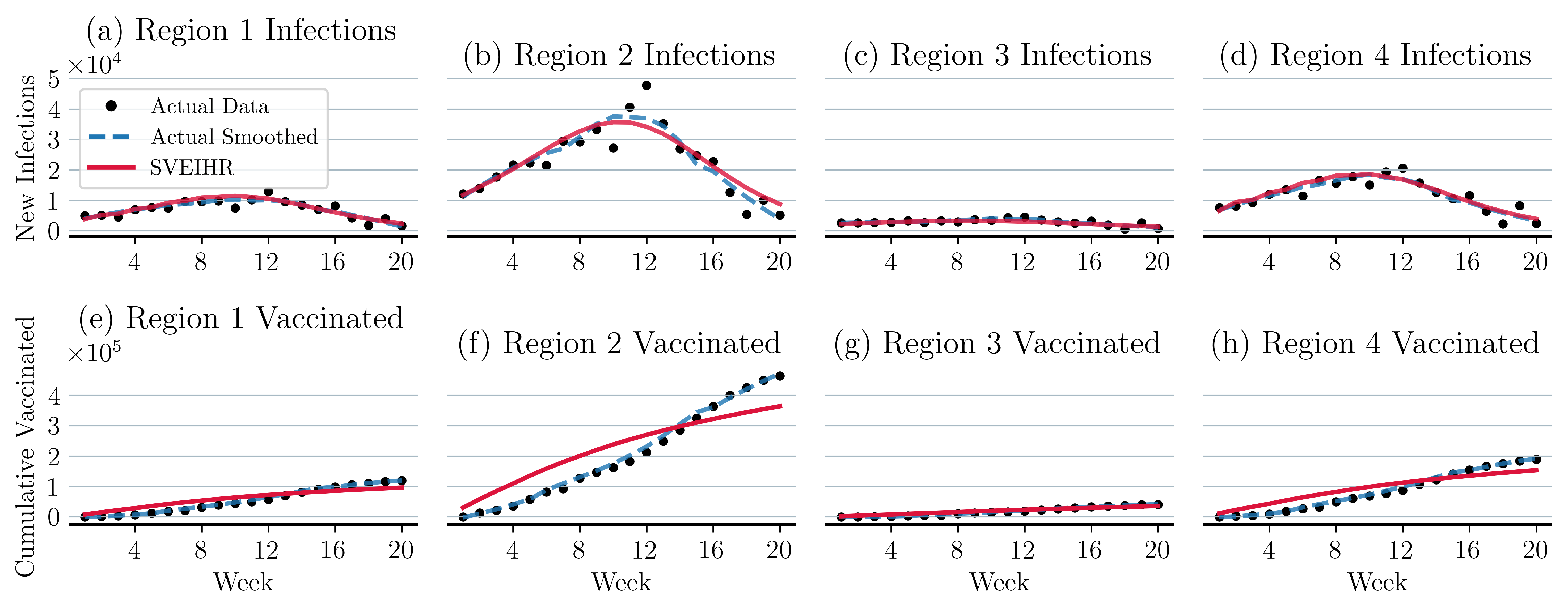}
    \caption{Comparison of predictions from the calibrated SVEIHR model and actual data across four regions for the baseline scenario, illustrating new infections and cumulative vaccinations over time.}
    \label{fig:calibration}
\end{figure}

Figure~\ref{fig:calibration} summarizes the results of our model calibration of the baseline scenario. Each plot in this figure presents the actual data and \ac{SVEIHR} model predictions of the weekly number of infections and vaccination intakes for each region of study. New infections at time $t$ are calculated using this equation: $\alpha(p^{\text{m}}+p^{\text{s}})E_{t,r} + \alpha(p^{\text{m},v} +p^{\text{s},v})EV_{t,r}$. The cumulative number of vaccinations up to week $t$ is calculated using $\sum_t (1-h_{t,r})\rho_r S_{t,r}$. We employ the Optuna framework to optimize parameter tuning \citep{optuna}. Table~\ref{tab:calibrate} summarizes the different performance metrics that we use to evaluate the performance of \ac{SVEIHR} model of the baseline scenario, including \ac{RMSE}, \ac{MAE}, and \ac{MAPE}. The low value of \ac{MAPE} highlights the accuracy of \ac{SVEIHR} in tracking new infections and vaccination trends in each region.

\begin{table}
    \centering
    \caption{Statistical analysis of prediction error from the calibrated \ac{SVEIHR} model for the baseline scenario.}
    \label{tab:calibrate}
    \begin{tabular}{c|ccc|ccc}
    \toprule
    & \multicolumn{3}{c|}{\textbf{Infections: Fig.~\ref{fig:calibration}(a) - (d) }} & \multicolumn{3}{c}{\textbf{Vaccinations: Fig.~\ref{fig:calibration} (e)-(h)}} \\
    \textbf{Region} & \textbf{RMSE} & \textbf{MAE} & \textbf{MAPE} & \textbf{RMSE} & \textbf{MAE} & \textbf{MAPE}\\
    \midrule
    R1 & 705.43 & 545.28 & 0.10\% & 17,094.55 & 15,883.84 & 2.00\% \\
    R2 & 2,288.45 & 1,877.74 & 0.16\% & 62,287.10 & 57,541.79 & 1.50\% \\
    R3 & 394.94 & 292.13 & 0.10\% & 4,531.79 & 4,177.05 & 4.06\% \\
    R4 & 765.10 & 620.95 & 0.06\% & 26,119.66 & 24,296.23 & 2.11\% \\
    \bottomrule
    \end{tabular}
\end{table}

\subsection{\emph{Value of stochastic solution}} \label{s:discussion.vss}

The \acf{VSS} quantifies the benefit of incorporating uncertainty into decision-making processes by comparing the results of the compartmental MSP with those of its deterministic counterpart. To compute the VSS for each decision stage $j \in \mathcal{J}$, we adopt the methodology proposed by \cite{escudero_2007}. In this method, the \ac{EV} model, which operates on the expected values of uncertain parameters, is used as the deterministic benchmark. Let $\mathit{EEV}_j$ denote the optimal value of the proposed compartmental \ac{MSP} model when the solution up to the decision period immediately preceding $j$ is fixed to the solution from the deterministic model. Let $\mathcal{Z}^*$ denote the optimal objective function value of the compartmental \ac{MSP} model. Therefore, \ac{VSS} at stage $j$, represented as $\mathit{VSS}_j = \mathit{EEV}_j - \mathcal{Z}^*$, captures the value derived from accommodating uncertainty in the model up to decision stage $j$. 

\begin{figure}[htp]
    \centering
    \includegraphics[width=1\textwidth]{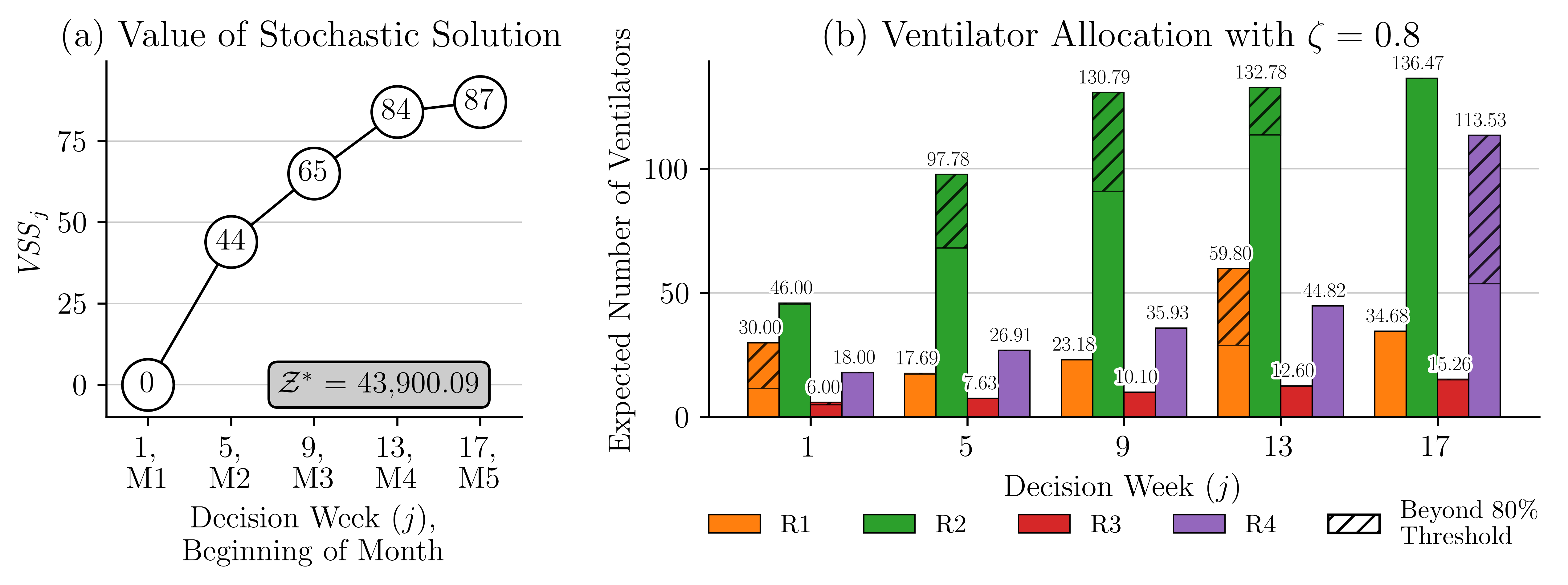}
    \caption{(a) Value of stochastic solution of the proposed model. (b) Expected ventilator allocation to different regions.}
    \label{fig:rq1}
    \footnotesize
    \floatfoot{\textit{Note.} Compartmental \ac{MSP} model include constraint~(\ref{eq:proportional}) with $\zeta$ is set to 0.8. The total number of ventilators available in each decision period: $\Delta_1=100, \Delta_5=150, \Delta_9=200, \Delta_{13}=250, \Delta_{17}=300$.}
\end{figure}

The evolution of the \ac{VSS} across the decision stages is illustrated in Figure~\ref{fig:rq1}(a). For the model used in this illustration, there are 5 decision stages that correspond to weeks 1, 5, 9, 13, and 17. The results of this figure highlight the benefits derived from the proposed compartmental \ac{MSP} model. The \ac{VSS} starts at zero, which is expected since $\mathit{EEV}_1 = \mathcal{Z}^*$. As the model progresses through the stages, the \ac{VSS} trends upward, reflecting the benefits of incorporating uncertainty into strategic decision making. We notice a decrease in the rate of change in \ac{VSS} from week 13 to 17. This is mainly because the mean change in VH is lower during these weeks. The results in Figure~\ref{fig:modeling_vh} show that the mean change in VH from week 13 to week 17 ranged from 0.03 to 0.05, which is lower than 0.06-0.09 seen from weeks 9 to 13.

Note that $VSS_j = 87$ in the illustrative example in Figure~\ref{fig:rq1}(a) appears to be a rather small number compared to the optimal solution $\mathcal{Z}^* = 43,900.$  In this example, a critically ill patient uses a ventilator for three weeks, the recovery rate is 51\% and $\overline{p} = 300$. This leads to a total of 1,000 ventilators available for allocation and use over five months, and a maximum of 729 critically ill patients to be saved using these ventilators. Thus, a $VSS_j = 87$ represents 12.21\% of this total, which is pretty high and highlights the value of using the compartmental MSP.  

\section{Discussion of the Results} \label{s:discussion}

\noindent{{\bf RQ1:} \emph{Can we develop an effective resource allocation model of critical resources needed during disease outbreaks?}} \label{s:discussion.rq1}

To respond to this research question, we develop the compartmental MSP model that facilitates resource allocation decisions under uncertainty. We make three observations: ($i$) \emph{The compartmental MSP model facilitates capacity allocation decisions given the progression and spread of a disease.} One could also simulate the SVEIHR model for different values of the number of ventilators available during the planning period and for different assignments of these ventilators to regions of study. However, such an approach would be time-consuming. ($ii$) \emph{The MSP model, rather than its deterministic counterpart, leads to better decision making under uncertainty.} We demonstrate this via our discussion of the VSS in Section~\ref{s:discussion.vss}.  ($iii$) \emph{The MSP model facilitates incorporating fairness in capacity allocation decisions.} 

Let us demonstrate these observations using the example in Figure~\ref{fig:rq1}(b). The figure shows the allocation of ventilators among our four study regions over a 20-week planning period. The decision periods are weeks 1, 5, 9, 13 and 17. In this experiment, we use an equity threshold $\zeta = 0.8$, which enforces that 80\% of the ventilators are allocated based on the population size of each region. The remaining 20\% of the ventilators are allocated in a way that minimizes the objective function. To help visualize the results, we use a different color for each region. The hatched portions of the bars represent the allocation of the remaining 20\% of the ventilators. These results indicate that in the early stages (weeks 5 and 9), the model allocates the additional ventilators to region 2 since this region's population size and density are the highest. Its population corresponds to 56.84\% of the total. Allocating additional ventilators to this region has the greatest impact on reducing the number of deaths. Later in the planning period, additional ventilators are allocated to regions 1 and 4. This is due to an increase in COVID-19 cases and increased vulnerability within these regions. These dynamic adjustments of resource allocation emphasize the model's ability to respond to changing needs throughout the pandemic,  providing an efficient and fair distribution of resources. 

These observations highlight the value that the \ac{MSP} model brings to decision makers during health crises, enabling them to balance between efficiently and equitably allocating scarce critical resources.

\noindent{{\bf RQ2:} \emph{What should be effective resource allocation strategies at the start of an outbreak?}} \label{s:discussion.rq2}

\begin{figure}[htp]
    \centering
    \includegraphics[width=1\textwidth]{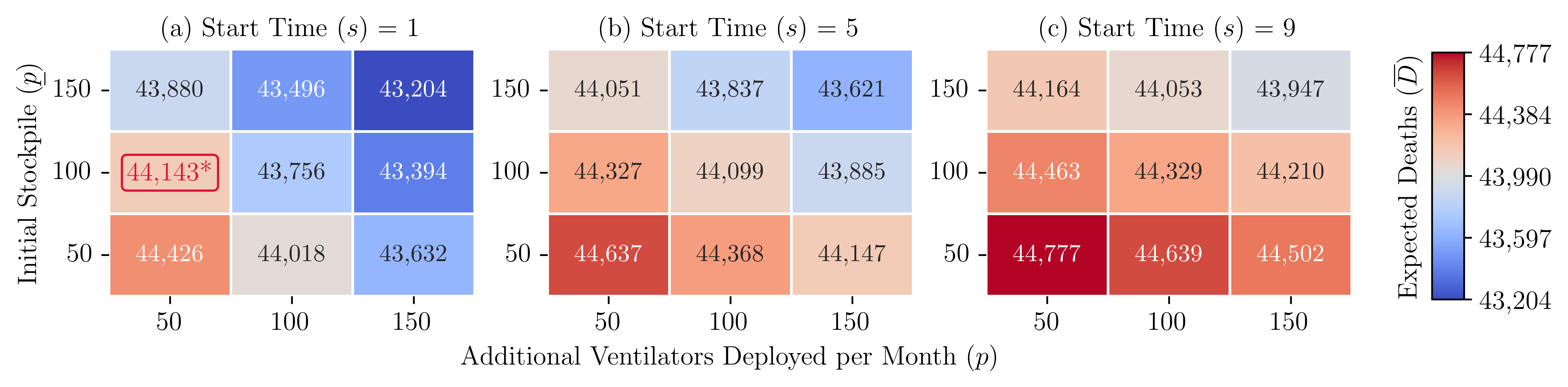}
    \caption{Impact of varying timing of the deployment and the number of new ventilators on the expected death toll.}
    \label{fig:rq2}
    \footnotesize
    \floatfoot{\textit{Note.} $^*$ The ventilator allocation is shown in Fig.~\ref{fig:rq1}.}
\end{figure}

To answer this research question, we designed the following experiment. We create three sets of problems. In each set, the time at which additional ventilators are deployed is different. In set 1, additional ventilators are deployed beginning in period 1. In set 2, that happens in period 5, and in set 3, it happens in period 9 of a 20-week planning period.  Each set consists of 9 problems, one for different combinations of the initial stockpile ($\underline{p}$) and the number by which we increase the number of ventilators deployed in each period ($p_j$). These experiments help us understand the impact of the timing and scope of decisions related to increasing the availability of critical resources (such as the timing and scope of the DPA) on the total expected number of deaths. A summary of the results from these experiments is presented in Figure~\ref{fig:rq2}. Detailed solution results for a problem from set 1 with $\underline{p}=100, p_j = 50$, are presented in Figure~\ref{fig:rq1}(b).

We make the following observations: ($i$) \emph{The timing of decisions related to the deployment of critical resources greatly impacts the expected number of deaths.} Based on our experiments, beginning the deployment of ventilators in the 1$^\text{st}$ rather than in the $9^\text{th}$ week of the planning period could lead to an 1.30\%  average of decrease in the expected number of deaths, which corresponds to 571 individuals. ($ii$) \emph{The size of the initial stockpile and the number by which we increase the availability of our critical resource every period greatly impact the expected number of deaths.} We develop a simple linear regression model to explore the relationships between the time of deployment of additional ventilators (${s}$), initial stockpile ($\underline{p}$), the number of additional ventilators deployed ($p_j$), and expected deaths ($\overline{D}$).  Since in our experiments, we keep $p_j$ constant over time, below, we use $p$ instead. The model is described by:

\begin{equation}
    \overline{D} = 44527.66 + 285.41 \; s - 5.44 \; \underline{p} - 4.81 \; p  \label{regressionModel}
\end{equation}

This regression model has an adjusted $R^2=0.9477$, which shows that the timing of ventilator deployment, initial stockpile, and the amount by which we increase resource availability explain approximately 94.77\% of the variation in the expected number of deaths. Based on our experiments, delaying the start of the deployment of additional ventilators by one month could lead to an average increase in the expected number of deaths by 285.41/month, highlighting the importance of prompt action. Each additional ventilator in the initial stockpile and each additional ventilator deployed lead to a decrease in the expected number of deaths by 5.44 (equivalent to 1.09/month) and 4.81 (0.962/month), respectively, highlighting the importance of maintaining a large stockpile and scalable production responses. \\

\noindent{{\bf RQ3:} \emph{What are the trade-offs between equity, equality, and effectiveness of resource allocation at the onset of a pandemic?}} \label{s:discussion.rq3}

To answer this research question, we designed the following experiment. We create four problems, P1 to P4. The compartmental MSP model without constraints~\eqref{eq:equity}, ~\eqref{eq:proportional} and ~\eqref{eq:equal_alloc} determines a resource allocation strategy that minimizes the expected number of deaths. This is our first problem (P1), which maximizes efficiency. The compartmental MSP model with constraints~\eqref{eq:equity} focuses on the equitable distribution of resources. As we change the value of the tolerance level $k \in (0, 1]$ from 0 to 1, the equity requirements become less stringent. This is our second problem (P2). The compartmental MSP model with constraints~\eqref{eq:proportional} is our third problem (P3), which also focuses on the equitable distribution of resources. As we change the value of the tolerance level $\zeta \in (0,1)$ from 0 to 1,  the equity requirements become rigid. The compartmental MSP model with constraints~\eqref{eq:equal_alloc} is our fourth problem (P4), which focuses on the equal distribution of resources.  The decision periods and the availability of resources are similar to the experiment described in {\bf RQ1}.   

\begin{figure}[htp]
    \centering
    \includegraphics[width=1\textwidth]{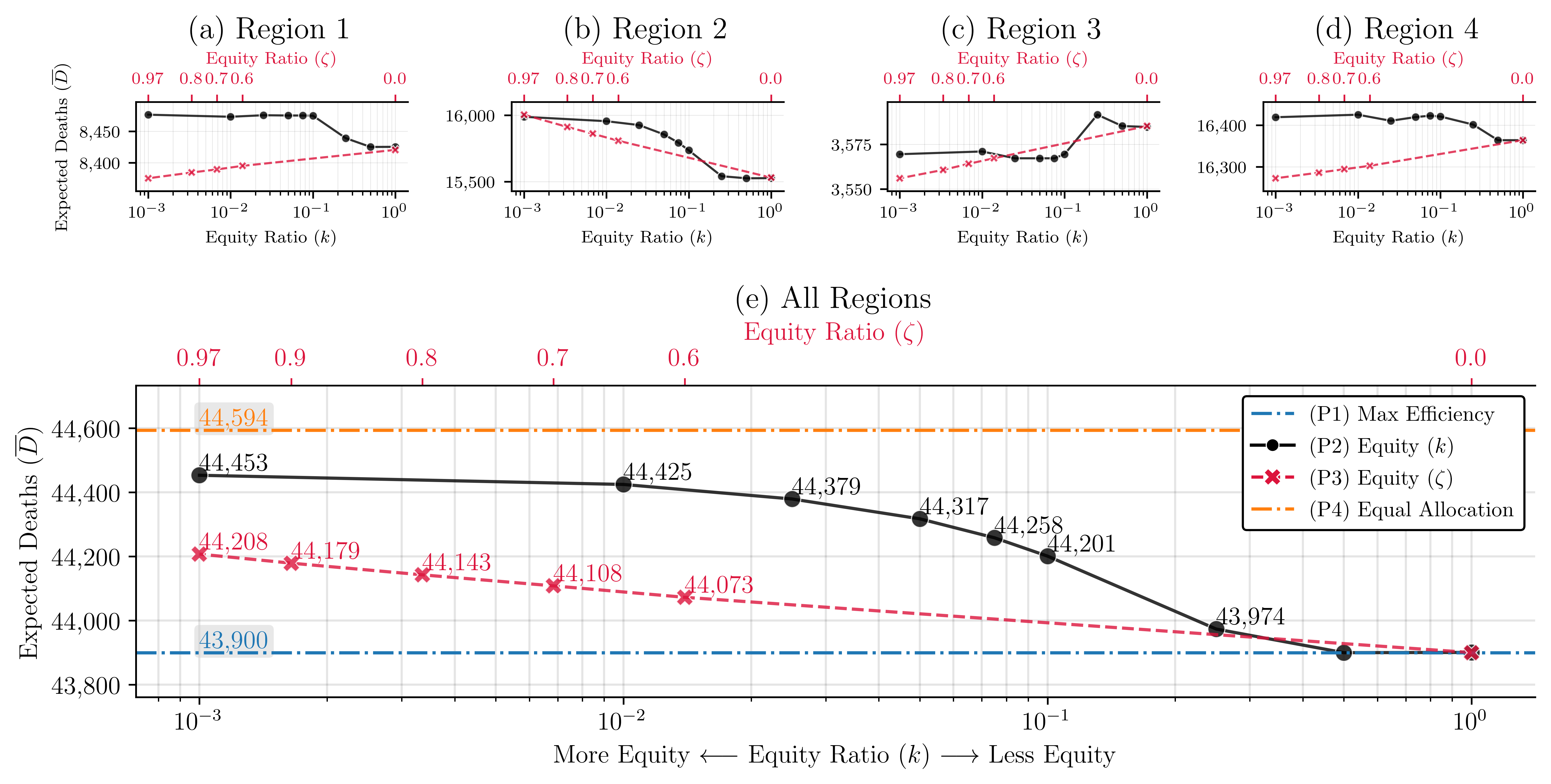}
    \caption{(a)-(d) Expected deaths at different equity levels for each region. (e) Total expected deaths for different resource allocation strategies.}
    \label{fig:rq3}
\end{figure}

We make the following observations: ($i$) \emph{Strategies that focus on equality underperform those that focus on efficiency and equity.}  Figure \ref{fig:rq3}(e) presents the total expected number of deaths as the value of $k$ changes from $10^{-3}$ to 1 (P2), and $\zeta$ changes from 0.97 to 0 (P3). Straight lines at 44,594 and 43,900 represent the expected deaths when ventilators are distributed equally (P4) and efficiently (P1), respectively. Based on these results, equal distribution of resources leads to the highest expected number of deaths. 

In our experiments, as the value of $k$ increases, the equity restrictions become less stringent, and the expected number of deaths decreases from 44,453 to 43,900, representing a 1.24\% decrease. Decreasing $\zeta$ makes equity restrictions less stringent, reducing the expected number of deaths from 44,208 to 43,900, which represents a 0.70\% decrease. The expected number of deaths is lowest when $k=1$, or $\zeta = 0$, which represents a model that focuses on efficiency (P1). 

($ii$) \emph{The price of ensuring equitable allocation of resources changes over time and space. This price is highest during the peak of a disease outbreak and in highly populated regions.} Figures~\ref{fig:rq3}(a) - (d) present the expected number of deaths in each region as the value of $k$ increases from $10^{-3}$ to 1, and the value of $\zeta$ decreases from 0.97 to 0. As $k$ increases and as $\zeta$ decreases, the equity restrictions become less stringent. These changes affect the regions of study differently. Region 3, which has the smallest population size and the lowest density, benefits the most when equity is enforced (the value of $k$ is close to $10^{-3}$). Region 2, which has the largest population size and the highest density, benefits the most when equity is not enforced (the value of $k$ is close to 1). Table \ref{tab:rq3} summarizes the expected number of deaths for each problem. The table also summarizes the change in the number of deaths for all problems compared to the results of P1.

\begin{table}
    \caption{Total expected deaths when applying different allocation strategies.}
    \label{tab:rq3}
    \begin{tabular}{ll@{\hskip 0.3in}cccccc}
    \toprule
    \multicolumn{2}{c}{\multirow{2}{*}{\textbf{Strategy}}}  & \multicolumn{6}{c}{\textbf{Expected Death}} \\
    \multicolumn{2}{l}{}  & \textbf{R1} & \textbf{R2} & \textbf{R3} & \textbf{R4} & \textbf{Total}  & \textbf{Diff from (P1)} \\
    \midrule
    \multicolumn{2}{l}{(P1) Max efficiency} & 8,427 & 15,524 & 3,585 & 16,364 & 43,900 & \_ \\
    \midrule
    \multirow{3}{*}{(P2) Equity ($k$)} & 0.001 & 8,477 & 15,988 & 3,569 & 16,419 & 44,453 & 553 \\
     & 0.05 & 8,475 & 15,855 & 3,567 & 16,420 & 44,317 & 417\\
     & 1 & 8,425 & 15,527 & 3,585 & 16,364 & 43,900 & 0 \\
    \midrule
    \multirow{3}{*}{(P3) Equity ($\zeta$)} & 0.97 & 8,375 & 16,004 & 3,556 & 16,272 & 44,208 & 308 \\
     & 0.8 & 8,384 & 15,912 & 3,561 & 16,286 & 44,143 & 243 \\
     & 0.6 & 8,395 & 15,808 & 3,567 & 16,303 & 44,073 & 173\\
    \midrule
    \multicolumn{2}{l}{(P4) Equal allocation} & 8,349 & 16,480 & 3,505 & 16,260 & 44,594 & 694 \\
    \bottomrule
    \end{tabular}
\end{table}

\begin{figure}[htp]
    \centering
    \includegraphics[width=0.7\textwidth]{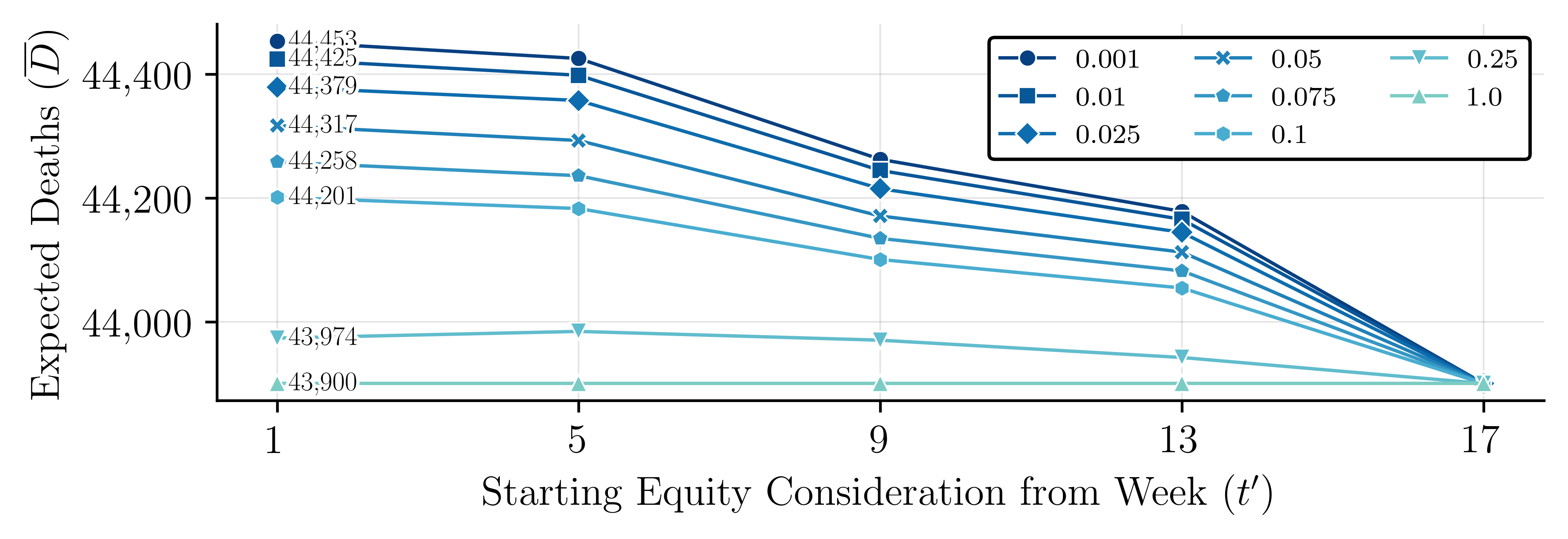}
    \caption{Evaluating the timing of resource allocation decisions on the expected number of deaths.}
    \label{fig:rq3_additional}
    \footnotesize
    \floatfoot{\textit{Note.} Summary of results from solving P2 at different time $t^{\prime}$.}
\end{figure}

To evaluate the impact that the timing of resource allocations has on the expected number of deaths, we performed the following experiment. We solve P2 several times by changing the time period over which we enforce the equity constraints~\eqref{eq:equity}. Thus, instead of adding for ${t\in\mathcal{T}}$, we add over ${t\in\{\mathcal{T}|t\geq t^{\prime}\}}$, which means that the equity in resource allocation is enforced after $t^{\prime}$. Prior to $t^{\prime}$, resource allocation is focused on maximizing efficiency.  We solve P2 for different values of $t^{\prime}$ and for different values of $k.$  The results of these experiments are summarized in Figure \ref{fig:rq3_additional}. These results show a decreasing trend in the expected number of deaths as $k$ approaches 1 and $t^{\prime}$ approaches $|\mathcal{T}|.$

This analysis is important for decision makers considering the appropriate level of equity tolerance that harmonizes both public health objectives and ethical standards. Our analysis provides a data-driven framework for such decisions, emphasizing the necessity of a comprehensive approach that aims to balance the imperative of saving lives and the principles of fairness and justice in resource allocation.

\section{Conclusion and future research directions} \label{s:conclusion}

\textbf{Summary of the Proposed Research:} This study proposes a modeling framework that integrates a compartmental model within a \acf{MSP} to tackle the complex challenge of allocating scarce healthcare resources during infectious disease outbreaks. The compartmental model simulates the dynamic changes of \acf{VH} over time and space, improving the ability of \ac{MSP} to adapt resource allocation strategies to the dynamics of disease spread and changes in population behavior toward vaccines. This compartmental MSP model also allowed us to analyze the price of fairness in resource allocation. Through our numerical analysis, we gain a deeper understanding of the impacts that the deployment time and size of the resource allocation, the effectiveness in managing the spread of the disease, and the price of fairness in resource allocation have on the expected number of deaths. The results of our analysis highlight the importance of advanced planning and preparedness for future outbreaks, effective management of critical resources, and the balance of equity and effectiveness in resource allocation.

\noindent\textbf{Research Contributions:} The main contribution of this work is the development of a compartmental MSP model that captures the impact of \ac{VH} in the effectiveness of allocating critical healthcare resources during a disease outbreak. The results of our numerical analysis show how public behaviors towards the COVID-19 vaccine impacted disease dynamics and the need for healthcare resources.  

Our review of the literature showed that there is very limited research on compartmental MSP models. Our proposed data-driven MSP compartmental model fills this gap in the literature and provides a robust framework for planning responses to disease outbreaks.  The model supports dynamic decision-making over time and space, taking into account the unpredictable nature of disease spread and resource availability. The model shows how dynamic, evidence-based strategies can enhance public health responses. The model addresses the price of fairness in the allocation of scarce healthcare resources. The model supports effective vaccination strategies and resource allocations that ensure greater immunity and mitigate the impact of outbreaks. 

\noindent\textbf{Research Findings:} Key observations from our study include: ($i$) The proposed compartmental \ac{MSP} model offers a more effective framework for resource allocation than an MSP model, as it captures the impact of disease progression (accounting for varying levels of \ac{VH}) on the need for critical healthcare resources. ($ii$) The proposed model supports resource allocation in different regions over time by considering both population size and immediate needs for medical resources, and balancing fairness and efficiency in decision making processes. ($iii$) The size of the initial stockpile, the timing of deployment of additional resources, and the number of additional resources have a significant impact on the expected number of deaths. Based on our numerical analysis, each unit increase in the initial stockpile and each unit increase of available ventilators could lead to a decrease in the expected number of deaths by 1.09/month and 0.962/month, respectively.  Based on our numerical analysis, delaying the deployment of critical resources by a month could increase the expected number of deaths by 285.41/month. ($iv$) An equal resource allocation strategy results in the highest number of expected deaths (that is, 44,594 over 5 months, assuming that 1,000 become available during this period), while a utilitarian resource allocation results in the smallest (i.e., 43,900 over 5 months). ($v$) The timing of the implementation of equity in resource allocation greatly impacts health outcomes. Delaying its implementation reduces deaths to levels comparable to those achieved via a utilitarian approach.

\noindent\textbf{Future Research Directions:} The model presented in this study can be extended in multiple ways to address other relevant problems related to the allocation of critical healthcare resources. The compartmental \ac{MSP} models could be extended to capture the interactions between various health interventions (i.e., masking mandates, travel restrictions, school closures, etc.) and the allocation of critical resources under uncertainty. These interventions affect the distribution of VH and the infection rate over time and space. Assume that the VH and infection rate are the uncertain problem parameters considered. In this case, the distribution of random variables is affected ``endogenously" by decisions made by the model. MSPs with endogenous uncertainty are more difficult to solve compared to models with exogenous uncertainty. 

Our proposed model can be extended to support other comparative studies across different regions or demographics to understand the effectiveness of strategies for allocating critical healthcare resources. In addition to the resources considered in our study, the proposed model could be extended to evaluate the impact of allocation strategies for other critical health resources (such as personal protective equipment, testing kits, and hospital beds) on the spread of the disease and the expected number of deaths. 

These model extensions provide opportunities to develop novel solution approaches and address relevant healthcare-related problems. 

\newpage

\begin{appendices}
\section{Impact of Healthcare Resources on Disease Outcomes}\label{appendix:compartmental_model}

\begin{figure}[htp]
    \centering
    \includegraphics[width=0.8\textwidth]{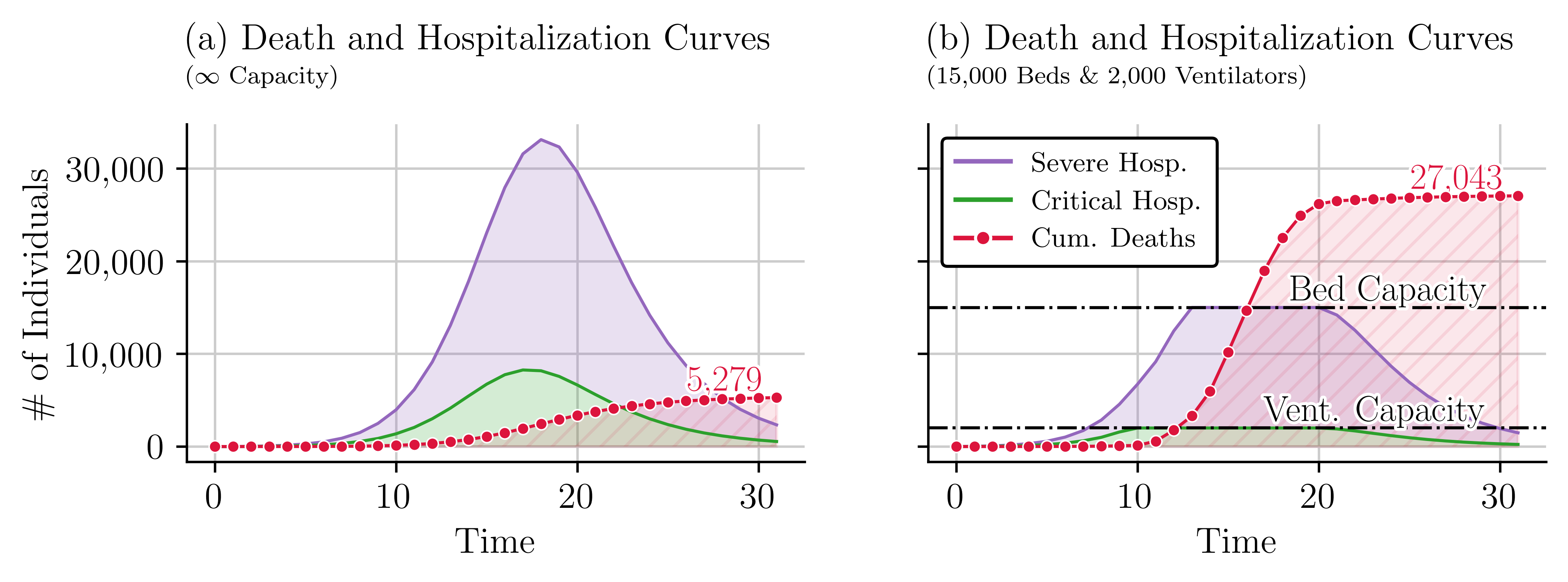}
    \caption{Comparative impact of healthcare resource availability on the SVEIHR model outcomes, illustrating differences in hospitalization and death rates under unlimited and limited resources scenarios.}
    \label{fig:sveihr_impact}
\end{figure}

The availability of healthcare resources such as beds and ventilators plays a critical role in managing the outcomes of infectious diseases. We performed two simulations of the \ac{SVEIHR} model to explore these relationships. The first simulation assumes unlimited healthcare resources, and the second introduces bed and ventilator availability constraints to reflect a real-world scenario with limited resource availability. The outcomes of the first simulation, depicted in Figure~\ref{fig:sveihr_impact}(a), show a low number of hospitalizations and deaths due to the assumption of infinite resources. However, the results of the second simulation, presented in Figure~\ref{fig:sveihr_impact}(b), reveal a higher number of deaths due to limited healthcare resources. The increase in the cumulative number of deaths becomes apparent when the number of individuals in critical condition exceeds the number of beds available. The contrast is notable. When the number of beds is unlimited, the expected number of deaths is 5,239, whereas when the number of beds is limited, the model results in 28,507 deaths.

\section{Modeling Uncertainty in Vaccine Hesitancy}\label{appendix:modeling_vh}

\begin{figure}[htp]
    \centering
    \includegraphics[width=\textwidth]{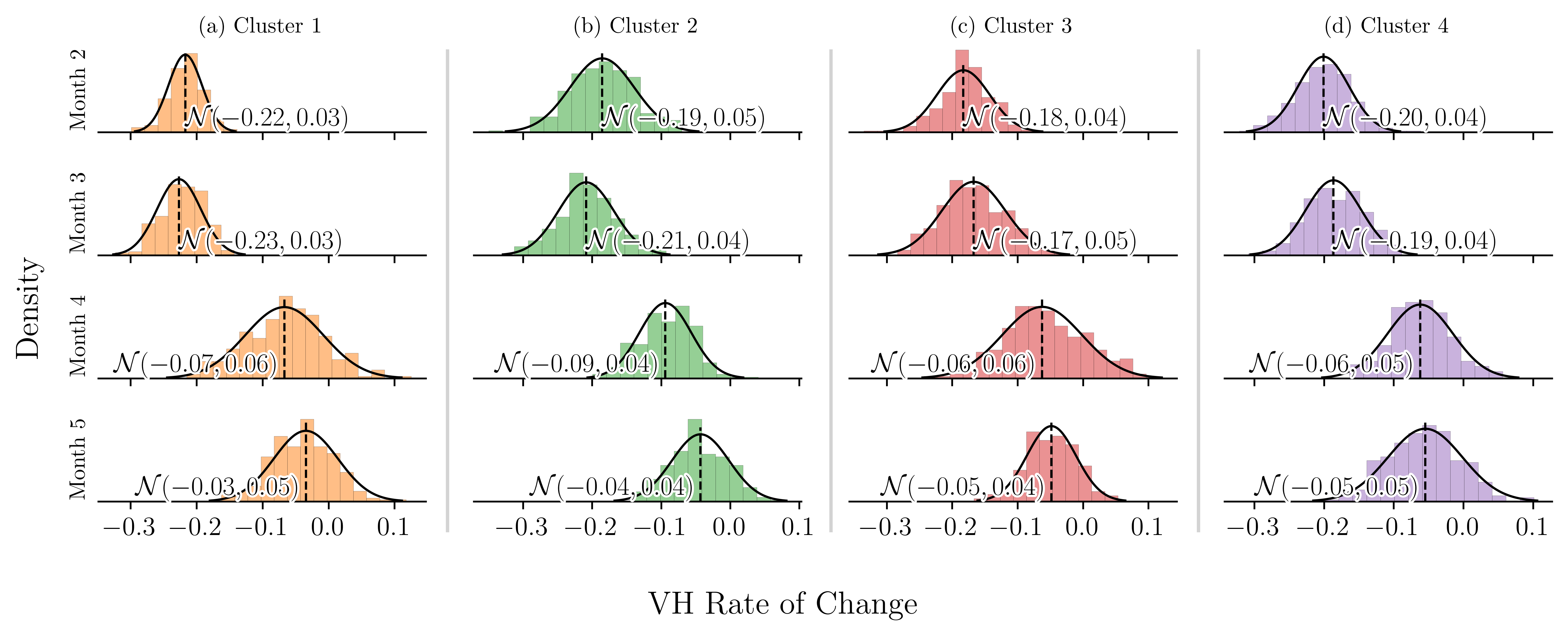}
    \caption{(a)-(d) Distribution of the monthly rate of change of vaccine hesitancy for different clusters in the study, with each row representing transitions from the previous month.}
    \label{fig:modeling_vh}
\end{figure}

Figure~\ref{fig:modeling_vh} shows that the rate of change in vaccine hesitancy fits the Normal Distribution. The corresponding mean and standard deviation changes over time. 

\section{Model Notation}\label{appendix:notation}

Tables~\ref{tab:sets}, \ref{tab:variables}, and \ref{tab:parameters} represent the model notations that are used throughout the paper.

\begin{table}
    \centering
    \caption{Sets and indices in the proposed \ac{MSP}.}
    \label{tab:sets}
    \begin{tabular}{ccl}
        \toprule
        \textbf{Set} & \textbf{Index} & \textbf{Description}\\
        \midrule
        $\mathcal{T}$ & $t$ & Set of time periods. \\
        $\mathcal{J}$ & $j$ & Set of decision time periods, $\mathcal{J}\subset\mathcal{T}$. \\
        $\mathcal{R}$ & $r$ & Set of regions. \\
        $\Omega$ & $\omega$ & Set of scenarios. \\
        $\mathcal{L}$ & $\ell$ & Set of state variables after fitting the SVEIHR model with actual data. \\
        $\mathcal{S}_{\omega,j,r}$ & & Set of scenarios that are indistinguishable with $\omega$ at decision stage $j$. \\
        \bottomrule
    \end{tabular}
\end{table}

\begin{table}
    \centering
    \caption{Variables in the proposed \ac{MSP}.}
    \label{tab:variables}
    \begin{threeparttable}
    \begin{tabular}{cl}
    \toprule
    \textbf{Variable} & \textbf{Description} (for $r\in\mathcal{R}$, at time $t \in \mathcal{T}$, $j\in\mathcal{J},$ and under scenario $\omega\in\Omega$)\\
    \midrule
    ${x_{\omega,j,r}}^{(1)}$  & Ventilators allocated at the start of decision period $j$. \\
    $\mathcal{X}_{\omega,t,r}$ & Cumulative numbers of ventilators at the end of time $t$. \\
    $S_{\omega,t,r},\; V_{\omega,t,r}$ & Number of susceptible individuals who are unvaccinated and vaccinated, respectively. \\
    $E_{\omega,t,r},\; EV_{\omega,t,r}$ & Number of exposed individuals who are unvaccinated and vaccinated, respectively. \\
    $I^{\text{m}}_{\omega,t,r},\; I^{\text{s}}_{\omega,t,r}$ & Number of infectious individuals with mild and severe symptoms, respectively. \\
    $H^{\text{s}}_{\omega,t,r},\; H^{\text{c}}_{\omega,t,r}$ & Number of hospitalized individuals with severe and critical conditions, respectively. \\
    $\mathcal{A}^{\text{s}}_{\omega,t,r},\; \mathcal{A}^{\text{c}}_{\omega,t,r}$ & Number of hospital admissions for severe and critical cases, respectively. \\
    $\mathcal{K}^{\text{s}}_{\omega,t,r},\; \mathcal{K}^{\text{c}}_{\omega,t,r}$ & Number of severe and critical cases not admitted to the hospital due to capacity constraints. \\
    $R_{\omega,t,r}$ & Cumulative number of recovered individuals up to period $t$.\\
    $D_{\omega,t,r}$ & Cumulative number of deaths up to period $t$.\\
    \bottomrule
    \end{tabular}
    \begin{tablenotes}\footnotesize
        \item[(1)] decision variable
    \end{tablenotes}
    \end{threeparttable}
\end{table}

\begin{table}
    \centering
    \caption{Parameters in the proposed \ac{MSP}.}
    \label{tab:parameters}
    \begin{tabular}{cl}
        \toprule
        \textbf{Parameter} & \textbf{Description} \\
        \midrule
        ${h_{\omega,t,r}}$& Vaccine hesitancy rate at region $r$ and time $t$ under scenario $\omega$. \\
        $\beta_r$ & COVID-19 contact rate in region $r$. $\theta_{\omega,t,r} = \beta_r S_{\omega,t,r}\left(I^{\text{m}}_{\omega,t,r} + I^{\text{s}}_{\omega,t,r}\right)/n_r$. \\
        $\rho_{r}$ & Maximum vaccination rate in region $r$. \\
        $\epsilon$ & Vaccine efficacy. \\
        $\alpha$ & Incubation rate. \\
        $\gamma^\text{m}_r,\; \gamma^\text{s}_r$ & Recovery rates for mild and severe cases in region $r$. \\
        $\gamma^\text{ks}_r,\; \gamma^\text{c}_r$ & Recovery rates for non-hospitalized severe and critical cases in region $r$. \\
        $\sigma_r$ & Rate from symptoms to hospitalization in region $r$. \\
        $\mu^{\text{ks}}_{r},\; \mu^{\text{c}}_{r}$ & Mortality rates for non-hospitalized severe and hospitalized critical cases in region $r$. \\
        $\varsigma^{\text{ks}}_r,\; \varsigma^{\text{c}}_r$ & Recovery rates for non-hospitalized severe and hospitalized critical cases in the region $r$. \\
        $k$ & Equity threshold in constraint~\eqref{eq:equity}. \\
        $\zeta$ & Equity threshold in constraint~\eqref{eq:proportional}.\\
        $p^{\text{m}},\; p^{\text{m},v},\; p^{\text{s}},\; p^{\text{s},v}$ & Transition percentages to mild and severe cases for unvaccinated and vaccinated individuals.  \\
        $p^{\text{r},v}$ & Recovery rate for vaccinated individuals. \\
        $\boldsymbol{\pi}_{l,r}$ & Initial count of individuals in state $\ell$ (compartment in \ac{SVEIHR} model) at time $t=0$ in region $r$.\\
        $n_r$ & Population size of region $r$. \\
        $b_r$ & Hospital bed count in region $r$. \\
        $\Delta_t$ & Additional ventilators available to allocate at start of time $t$. \\
        \bottomrule
    \end{tabular}
\end{table}

\section{Constraint Linearization}\label{appendix:linearize}

The proposed compartmental \ac{MSP} model presented in Section~\ref{s:method.compartmental_msp} contains three non-linear constraints. They can be linearized as follows:

We linearize the \textbf{`min'} functions in constraints~\eqref{eq:admit_critical} and \eqref{eq:admit_severe} using the big-M method. To illustrate, let's consider the linearization of constraint~\eqref{eq:admit_severe}. For ease of demonstration, let us declare a few auxiliary variables: $a_1=\sigma I^{\text{s}}_{\omega,t,r}$ representing the number of severe case admissions and $a_2= b_{r} - H^{\text{c}}_{\omega,t,r} - H^{\text{s}}_{\omega,t,r} - \mathcal{A}^{\text{c}}_{\omega,t,r}$ representing the remaining bed capacity after accounting for current hospitalizations. 

The original constraint can be written as $\mathcal{A}^{\text{s}}_{\omega,t,r} = \min\{a_1, a_2\}$, indicating that the number of severe cases admitted is the lesser of the demand or the remaining capacity. We use a binary variable ($y_{\omega,t,r}$) and big-M constant ($M$) to enforce $A\geq \min\{a_1, a_2\}$ . The original \textbf{`min'} constraint is then replaced by the following system of linear inequalities~\eqref{eq:linear_min_1}-\eqref{eq:linear_min_2}:

\begin{subequations}
    \label{eq:linear_min}
    \renewcommand{\theequation}{\ref{eq:admit_severe}\alph{equation}}
    \begin{align}
        \mathcal{A}^{\text{s}}_{\omega,t,r} &\leq a_1 \label{eq:linear_min_1}\\
        \mathcal{A}^{\text{s}}_{\omega,t,r} &\leq a_2\\
        \mathcal{A}^{\text{s}}_{\omega,t,r} &\geq a_1 - M(1-y_{\omega,t,r})\\
        \mathcal{A}^{\text{s}}_{\omega,t,r} &\geq a_2 - My_{\omega,t,r} \label{eq:linear_min_2}
    \end{align}
\end{subequations}

Next, we linearize the absolute function in constraint~\eqref{eq:equity} using the following procedure. Again, for demonstration purposes, let us introduce the following auxiliary variables: $a_1=\sum_{t \in \mathcal{T}}\sum_{\omega\in\Omega}p_\omega H^{\text{c}}_{\omega,t,r}$, and $a_2=\sum_{t \in \mathcal{T}}\sum_{r'\in \mathcal{R}}\sum_{\omega\in\Omega} p_\omega H^{\text{c}}_{\omega,t,r'}$. Since $a_2 > 0, \sum_{r'\in \mathcal{R}}n_{r'} > 0,$ the original constraint can be written as $\left\lvert\frac{a_1}{a_2} - \frac{n_r}{\sum_{r'\in \mathcal{R}}N_{r'}} \right\rvert \leq k$ and is equivalent to:

\begin{align}
    &\begin{dcases*}
        \frac{a_1}{a_2} - \frac{n_r}{\sum_{r'\in \mathcal{R}}n_{r'}} \leq k\nonumber\\
        \frac{a_1}{a_2} - \frac{n_r}{\sum_{r'\in \mathcal{R}}n_{r'}} \geq -k\nonumber
    \end{dcases*}\\
    \Leftrightarrow &\begin{cases}
        a_1\left(\sum_{r'\in \mathcal{R}}n_{r'}\right) - a_2n_r \leq ka_2\left(\sum_{r'\in \mathcal{R}}n_{r'}\right)\nonumber\\
        a_1\left(\sum_{r'\in \mathcal{R}}n_{r'}\right) - a_2n_r \geq -ka_2\left(\sum_{r'\in \mathcal{R}}n_{r'}\right)\nonumber
    \end{cases}
\end{align}

Constraint~\eqref{eq:equity} can be replaced by the following linear inequalities~\eqref{eq:linear_abs_1}-\eqref{eq:linear_abs_2}:

\begin{subequations}
    \label{eq:linear_abs}
    \renewcommand{\theequation}{\ref{eq:equity}\alph{equation}}
    \begin{align}
    a_1\left(\sum_{r'\in \mathcal{R}}N_{r'}\right) - a_2n_r - ka_2\left(\sum_{r'\in \mathcal{R}}n_{r'}\right) &\leq 0\label{eq:linear_abs_1}\\
    a_1\left(\sum_{r'\in \mathcal{R}}N_{r'}\right) - a_2n_r + ka_2\left(\sum_{r'\in \mathcal{R}}n_{r'}\right) &\geq 0\label{eq:linear_abs_2}
    \end{align}
\end{subequations}

\section{Illustration of \acfp{NAC}}\label{appendix:nac}

\begin{figure}[htp]
    \centering
    \includegraphics[width=0.9\textwidth]{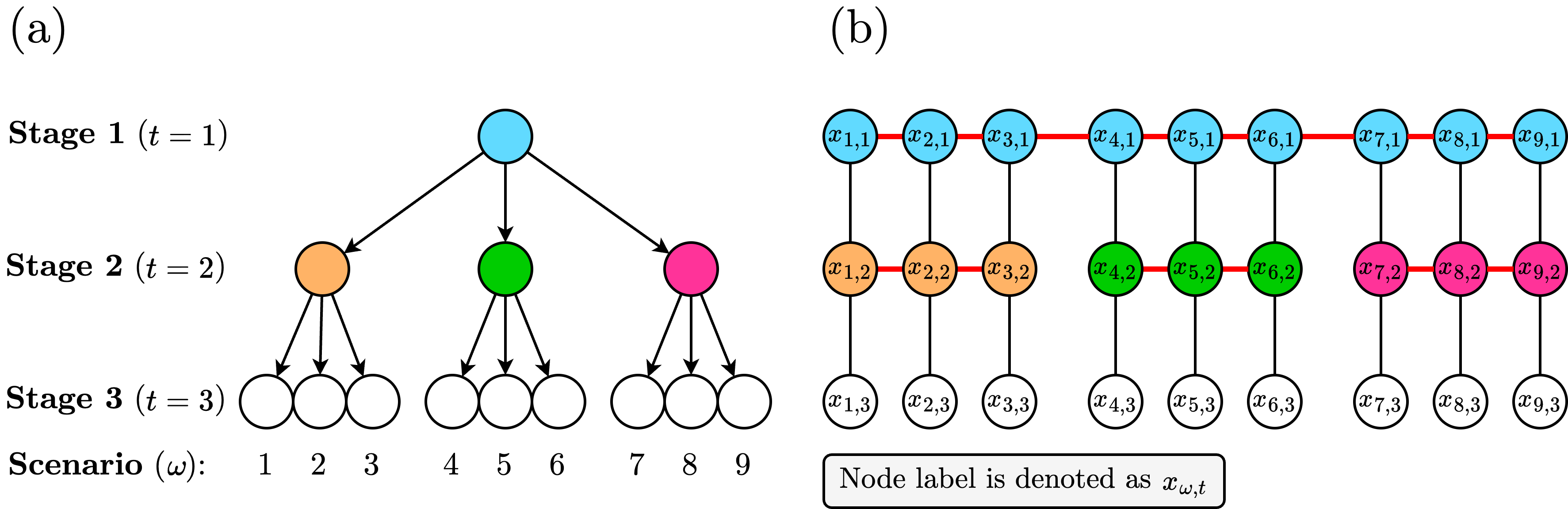}
    \caption{(a) Illustration of a scenario tree with three stages and nine scenarios, (b) an alternative form to visualize the \acp{NAC}, which are shown in red.}
    \label{fig:nac}
\end{figure}

Figure~\ref{fig:nac}(a) shows a scenario tree with three stages and nine possible scenarios. Figures~\ref{fig:nac}(a) and (b) help the reader understand the implications of the nonanticipativity constraints (NAC). \acp{NAC} ensure that decisions made at a certain stage are the same across all scenarios that share the same history up to that stage. For example, the root node (stage 1) in Figure~\ref{fig:nac}(a) is colored cyan. This node is common in all scenarios. As such, decisions made in stage 1 affect decisions made along all scenario paths. This is why in Figure~\ref{fig:nac}(b), the nodes in the top layer are all colored cyan. Therefore, the set of scenarios that are indistinguishable from $\omega = 1$ at time $t=1$ is $\mathcal{S}_{1,1,r} = \{2, 3, \cdots, 9\}$. We use the following eight NACs to ensure this: $x_{1,1,r} = x_{2,1,r};\; x_{1,1,r} = x_{3,1,r}; \dots;\; x_{1,1,r} = x_{9,1,r}$.

\section{Population Flow}\label{appendix:migration_flow}

\begin{figure}[htp]
    \centering
    \includegraphics[width=0.7\textwidth]{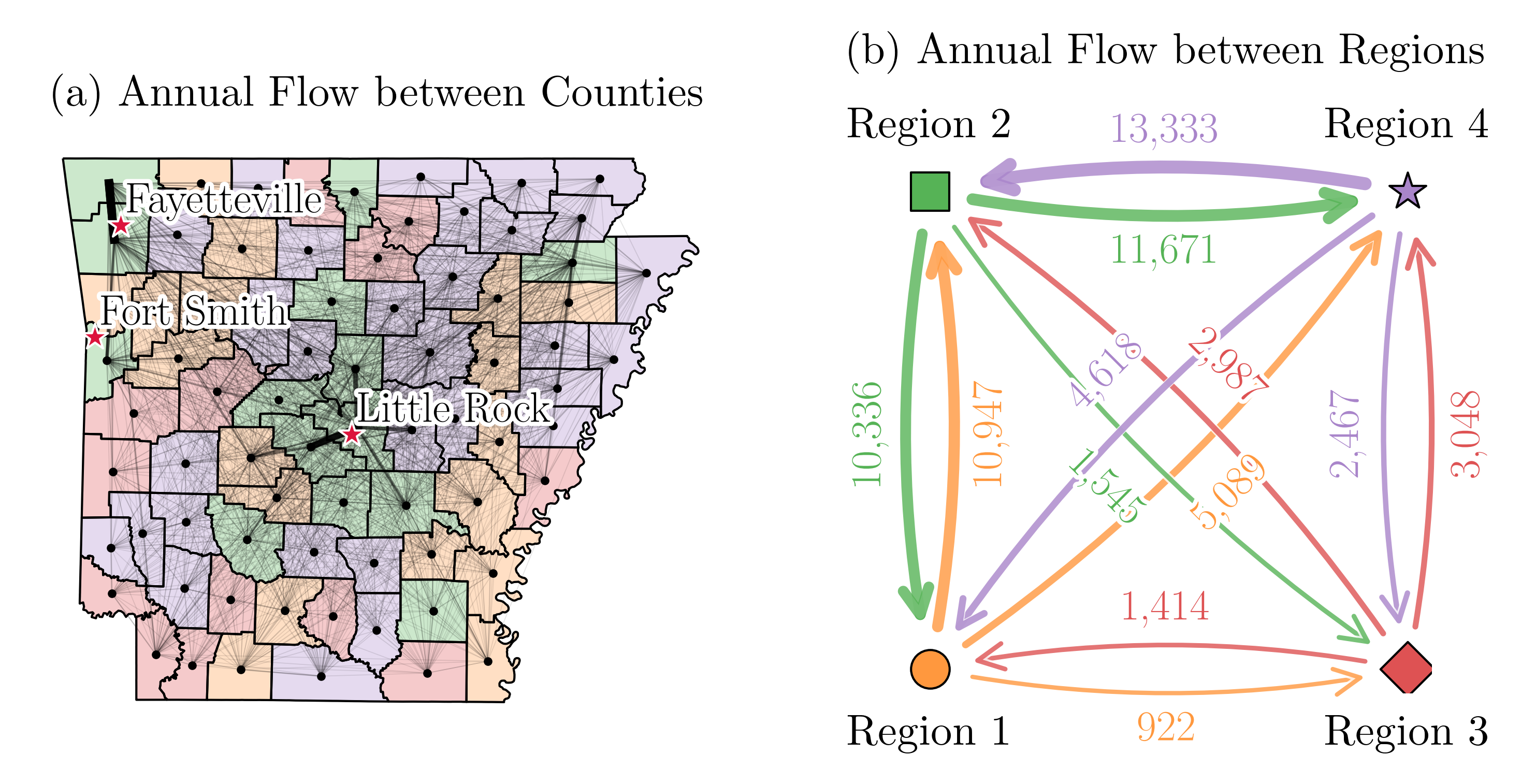}
    \caption{County and regional migration patterns in Arkansas: Analyzing annual inter-county flows and aggregate regional migration trends.}
    \label{fig:migration_flow}
\end{figure}

We utilized data related to county-to-county migration flows from 2016 to 2020, sourced from \citep{us_census_bureau}. We estimate the annual movement of individuals between counties in Arkansas. The results are shown in Figure~\ref{fig:migration_flow}(a). The black dots in this figure represent county centroids and the shaded areas represent the four different regions of study. Major cities like Fayetteville, Fort Smith, and Little Rock are labeled, signifying possible hubs of migration. The simplified flow between regions is depicted in Fig.~\ref{fig:migration_flow}(b). These data are integrated into the compartmental MSP model to update the susceptible individuals at each time step $t$ and region $r$ in constraint~\eqref{eq:susceptible}. The weekly movement rate are shown in Table~\ref{tab:migration_flow}. The expression $\sum_{r'\in \mathcal{R}} \nu_{r'\rightarrow r} S_{\omega,t,r'} - \sum_{r'\in \mathcal{R}} \nu_{r\rightarrow r'} S_{\omega,t,r}$ represents the net movement of susceptible individuals between regions. Since most migration occurs within regions, the flow between regions has a negligible impact on the model.

\begin{table}
    \centering
    \caption{Weekly population movement rate between studied regions.}
    \label{tab:migration_flow}
    \begin{tabular}{ccccc}
    \toprule
    \textbf{From \textbackslash To} & \textbf{R1} & \textbf{R2} & \textbf{R3} & \textbf{R4} \\
    \midrule
    \textbf{R1} & - & 5.34e-04 & 4.50e-05 & 2.48e-04 \\
    \textbf{R2} & 1.28e-04 & - & 1.92e-05 & 1.45e-04 \\
    \textbf{R3} & 1.58e-04 & 3.34e-04 & - & 3.41e-04 \\
    \textbf{R4} & 1.45e-04 & 4.19e-04 & 7.76e-05 & - \\
    \bottomrule
    \end{tabular}
\end{table}

\section{Model Data}\label{appendix:model_data}

Table~\ref{tab:model_data} presents the parameter values used in the proposed compartmental \ac{MSP} model. Table~\ref{tab:model_data_region} details the parameter values for the four different regions included in the study. The upper portion of the table includes the tuned parameters derived from the calibration process, which is described in Section~\ref{s:model.calibration}.  

\begin{table}
    \centering
    \caption{Parameter values used in the proposed model.}
    \label{tab:model_data}
    \begin{tabular}{ccc|ccc}
    \hline
    \textbf{Parameter} & \textbf{Value} & \textbf{Source} & \textbf{Parameter} & \textbf{Value} & \textbf{Source} \\
    \hline
    $|\mathcal{T}|$ & 20 &  & $\gamma^\text{s}, \gamma^\text{ks}, \gamma^\text{c}$ & 0.637, 0.233, 0.333 & \citep{Zhou2020} \\
    $|\mathcal{R}|$ & 4 &  & $\mu^\text{ks}$ & 1.000 &  \\
    $|\Omega|$ & 81 &  & $p^{\text{m}}, p^{\text{s}}, p^{\text{c}}$ & 0.800, 0.150, 0.050 & \citep{wu_mcgoogan_2020} \\
    $\mathcal{J}$ & $\{1,5,9,13,17\}$ & & $p^{\text{r},v}, p^{\text{m},v}, p^{\text{s},v}$ & 0.850, 0.100, 0.050 & \citep{Barnes2023} \\
    $\epsilon$ & 0.950 & \citep{kathy_2023} & $\varsigma^{\text{ks}}$ & 0.400 &  \\
    $\gamma^v$ & 1.000 &  & $\varsigma^{\text{c}}$ & 0.510 & \citep{wu_mcgoogan_2020}  \\
    \hline
    \end{tabular}
\end{table}

\begin{table}
    \centering
    \caption{Parameter values for each region used in the proposed model.}
    \label{tab:model_data_region}
    \begin{threeparttable}
    \begin{tabular}{ccccc}
    \toprule
    \multirow{2}{*}{\textbf{Parameter}} & \multicolumn{4}{c}{\textbf{Region ($\pmb{r}$)}}\\
     & \textbf{R1} & \textbf{R2} & \textbf{R3} & \textbf{R4} \\
    \midrule
    $\beta_r$ $^{(1)}$  & 1.873 & 2.162 & 1.805 & 1.781 \\
    $\rho_r$ $^{(1)}$& 0.051 & 0.051 & 0.043 & 0.050 \\
    $\gamma^{\text{m}}_r$ $^{(1)}$& 0.704 & 0.704 & 0.704 & 0.699 \\
    $\sigma_r$ $^{(1)}$& 0.943 & 0.877 & 0.952 & 0.901 \\
    $S_{0,r}$ $^{(1)}$& 299,000 & 1,107,000 & 129,000 & 502,000 \\
    $E_{0,r}$ $^{(1)}$& 3,100 & 3,200 & 100 & 3,500 \\
    $I^{\text{m}}_{0,r}$ $^{(1)}$& 200 & 5,900 & 1,300 & 3,600 \\
    $I^{\text{s}}_{0,r}$ $^{(1)}$& 2,500 & 1,470 & 80 & 60 \\
    \midrule
    $b_r$ & 1,405 & 9,357 & 490 & 2,032 \\
    $\mathcal{X}_{0,r}$ & 77 & 566 & 18 & 105 \\
    ${n_r}/{\sum n_r}$ & 0.145 & 0.568 & 0.063 & 0.224 \\
    \bottomrule
    \end{tabular}
    \begin{tablenotes}\footnotesize
        \item[(1)] Tuned parameters using Optuna framework, an open-source optimization library that automates hyperparameter tuning.
    \end{tablenotes}
    \end{threeparttable}
\end{table}

\end{appendices}

\bibliographystyle{unsrtnat}
\bibliography{ref}

\begin{acronym}
    \acro{OR}{operations research}
    \acro{SIR}{susceptible-infectious-recovered}
    \acro{SEIR}{susceptible-exposed-infectious-recovered}
    \acro{SVEIHR}{susceptible-vaccinated-exposed-infectious-hospitalized-recovered}
    \acro{SP}{stochastic programming}
    \acro{PPE}{personal protective equipment}
    \acro{ICU}{intensive care unit}
    \acro{WHO}{World Health Organization}
    \acro{MSP}{multi-stage stochastic program}
    \acro{VH}{vaccine hesitancy}
    \acro{CTIS}{COVID-19 Trends and Impact Survey}
    \acro{ODE}{ordinary differential equation}
    \acro{DTW}{dynamic time warping}
    \acro{CDC}{Centers for Disease Control and Prevention}
    \acro{NAC}{non-anticipativity constraint}
    \acro{RMSE}{root-mean-square error}
    \acro{MAPE}{mean absolute percentage error}
    \acro{MAE}{mean absolute error}
    \acro{GM}{General Motors}
    \acro{DPA}{Defense Production Act}
    \acro{HHS}{U.S. Department of Health and Human Services}
    \acro{VSS}{value of stochastic solution}
    \acro{EV}{expected value}
\end{acronym}

\end{document}